\documentclass[12pt]{article}
\usepackage{amssymb,amsmath,amscd}

\newtheorem{proposition}{Proposition}
\newtheorem{theorem}[proposition]{Theorem}
\newtheorem{lemma}[proposition]{Lemma}
\newtheorem{remark}[proposition]{Remark}

\def\sk{\vskip .4cm}
\def\noi{\noindent}
\def\om{\omega}

\def\al{\alpha}

\def\be{\beta}
\def\ga{\gamma}

\def\epsi{\varepsilon}

\def\th{\theta}

\def\pa{\partial}
\def\part{\partial}

\def\Rmia#1#2{ (R_{q,r})^{#1}_{~~~#2} }
\def\PA#1#2{ P^{#1}_{A~~#2} }

\def\Rh{{\hat R}}

\def\Rhat#1#2{ \Rh^{#1}_{~~~#2} }

\def\Rhatinv#1#2{ (\Rh^{-1})^{#1}_{~~~#2} }

\def\Rha{{\hat R}}

\def\T#1#2{ T^{#1}_{~~#2} }

\def\rminus{r^{-1}}

\def\D{\Delta}

\def\n2{{{N+1} \over 2}}
\def\ap{a^{\prime}}
\def\bp{b^{\prime}}

\def\sma#1{\mbox{\footnotesize #1}}

\def\Q.E.D.{\rightline{$\Box$}}

\def\Co{\mathbb{C}} % numeri complessi

\def\fidi{\hskip5pt \vrule height4pt width4pt depth0pt \par}

\def\id{{\rm id}}

\def\Re{{\mathbb{R}}}  % numeri reali
\def\S{\mathbb{S}}
\def\Sq{\mathbb{S}_q}
\def\Sqr{\mathbb{S}_{q,r}}

\def\sphere{\mathbb{S}}
\def\sph4{\Sigma_q^4}
\def\Tr{{\rm Tr}}

\def\u4{{\bf u}(4)}
\def\V{{\cal V}}

\def\*{\star}

\begin{document}

\title{
\vskip -0.90cm
%\vskip -2.10cm
%~~~~~~~~~~~~~~~~~~~~~~~~~~~~~~~~~~~~~~~~~~~~~~~~~~~
%{\small{LMU-TPW-01-09}} \\[2em]
On the Noncommutative Geometry \\ of Twisted Spheres}

\author{Paolo Aschieri$^1$, Francesco Bonechi$^{2,3}$}

\maketitle
{\vskip -.45cm}
\centerline{{ \small
{ ${ }^1$ Sektion Physik der Ludwig-Maximilians-Universit\"at}}}
\centerline{{\small Theresienstr. 37, D-80333 M\"unchen, Germany}}
\vskip .2cm
\centerline{{\small  ${ }^2$ INFN
Sezione di Firenze}} \centerline{{\small ${ }^3$ Dipartimento di
Fisica, Universit\`a di Firenze, Italy. }} \centerline{{\small
e-mail: aschieri@theorie.physik.uni-muenchen.de, bonechi@fi.infn.it.}}
\medskip
\begin{abstract}
We describe noncommutative geometric aspects of twisted
deformations, in particular of the spheres
in Connes and Landi \cite{ConnesLandi} and in
Connes and Dubois ${}^{\!}$Violette \cite{ConnesDuboisViolette},
by using the differential and integral calculus on these spaces
that is covariant under the action of their corresponding quantum
symmetry groups.
%$4$--dimensional case and then in
%of spheres and show that they coincide with
%$\theta$--deformations introduced for
%generic dimension.
We start from multiparametric deformations of
the orthogonal groups and related planes and spheres. We show that
only in the twisted limit of these multiparametric deformations
the covariant calculus on the plane gives by a quotient procedure
a meaningful calculus on the sphere. In this calculus the external
algebra has the same dimension of the classical one. We develop
the Haar functional on spheres and use it to define an integral on
forms. In the twisted limit (differently from the general
multiparametric case) the Haar functional is a trace and we thus
obtain a {\it cycle} on the algebra. Moreover we explicitely construct
the $*$--Hodge operator on the space of forms on the plane and then by
quotient on the sphere. We apply our results to even spheres
and we compute the Chern--Connes pairing between the character of
this cycle, {\it i.e.} a cyclic $2n$--cocycle, and the instanton
projector defined in \cite{ConnesDuboisViolette}.
\end{abstract}

%%%%%%%%%%%%%%%%%%%%%% 1 %%%%%%%%%%%%%%%%%%

\thispagestyle{empty}

\section{Introduction}
Noncommutative geometry is an active research field both in
Physics and Mathematics. There are many examples of
noncommutative spaces and many different philosophies behind them.
It is therefore instructive to study
the geometry of specific examples in order to gain a better
understanding of the relations between different approaches.

Quantum spaces related to {\it standard} quantum groups, {\it
i.e.} FRT (\cite{FRT}) matrix quantum groups or Drinfeld--Jimbo quantum
enveloping algebras, have a rich noncommutative structure and
display remarkable properties for many
respects (their connection with integrable models and invariants of knots, semiclassical structures, etc...),
but, on the other hand, present some triviality when considered under
the light of the noncommutative geometry \`a la Connes.
For instance, the $C^*\!$--algebra of the standard
Podle\`s quantum $2$--spheres  and of the quantum $4$-sphere of
\cite{NOI} are isomorphic, {\it i.e.} the topology cannot distinguish
inbetween the two different classical dimensions. Indeed the study
of \cite{MNW} of Podle\`s $2$-spheres shows that
the Hochschild
dimension, that corresponds to the commutative notion of dimension,
is zero and periodic cyclic cohomology is generated by
$0$--traces. The drop of dimension appearing in this example seems
to be a generic property of these kind of quantum spaces.

In \cite{ConnesLandi} a deformation $\Sq^4$ of the $4$--sphere was
introduced with the property that the Hochschild dimension equals
the commutative dimension. In \cite{ConnesDuboisViolette} this
program was generalized to any dimension: a complete
classification of $3$--dimensional spherical manifolds was given,
while for generic dimension a family of deformed spheres, called
{\it $\theta$--deformed} spheres, was introduced. In general
\cite{ConnesDuboisViolette} a $\theta$--deformation $M_\theta$ of
a manifold $M$ equipped with a smooth action of the $n$--torus
$T^n$ is determined by defining $C^\infty(M_\theta)$ as the
invariant subalgebra (under the action of $T^n$) of the algebra
$C^\infty(M\times T^n_\theta)\equiv C^\infty(M)\widehat{\otimes}
C^\infty(T^n_\theta)$ of smooth functions on  $M\times
T^n_\theta$; where $T^n_\theta$ is the noncommutative $n$-torus.
This construction  allows for an easy definition of smooth
differential forms and of the $*$--Hodge operator. Although this
calculus is not connected with a quantum group action, it
coincides with the covariant calculus.

In \cite{Varilly} (see also \cite{Sitarz}) it was shown that
$\sphere_q^4$ is an homogeneous space of a {\it twisted}
deformation of $SO(5)$. Twisted quantum groups quantize
Poisson--Lie groups defined by classical $r$-matrices living in
the Cartan subalgebra of the semisimple Lie algebra. A useful
point of view in order to understand the different geometry of the
standard and the twisted case is given by {\it multiparametric
deformations}. Multiparametric quantum groups were defined by
introducing in the FRT construction \cite{FRT} the multiparametric
$R$--matrices \cite{Schirrmacher}. These quantum groups depend on
two sets of parameters and contain as limit structures both the
standard deformation and the twisted deformation. Homogeneous
spaces with respect to multiparametric quantum groups are
constructed in the same way as for standard deformations. It is
easy to verify that the $4$--sphere of \cite{ConnesLandi} appears
as the limit, to the twisted structure case, of the
multiparametric orthogonal $4$--sphere (the $4$--sphere covariant
under the action of the multiparametric quantum $SO(5)$ group);
similarly the $\theta$-deformed spaces and groups of
\cite{ConnesDuboisViolette} can be obtained as special limits to the
twisted case of the corresponding multiparametric quantum spaces and groups.

A basic tool in order to analize the geometry of (multiparametric)
quantum spaces is the covariant differential calculus. This calculus
is constructed from the requirement of covariance
under the action of the corresponding orthogonal quantum  groups
\cite{Wess:1991vh, Carow-Watamura:1991zp}.
It is closely related to the
bicovariant noncommutative calculus on quantum groups \cite{Woronowicz}
(see also \cite{JurcoWat} and \cite{AschieriCastellani}).
In the past years the bicovariant calculus has been studied and
classification results are known.
Nevertheless a relation with the
noncommutative geometry \`a la Connes has lacked; in particular with
the notion of {\it cycle} over an algebra (\cite{ConnesBook}). Cycles
are differential algebras provided with an integral, which is a
closed graded trace; they are in one--to--one correspondence with
cyclic cocycles and allow the computation of topological
invariants like Chern characters.

In this letter we show that quantum groups techniques, and in
particular the covariant differential calculus on quantum planes
and spheres, are useful tools in order to describe the
noncommutative geometry of twisted deformations. In Section 2 we
introduce multiparametric deformations of the orthogonal groups
and of the related planes and spheres. The multiparametric
covariant differential calculus on the plane is then studied, we
see that only in the twisted limit the calculus admits a euclidean
reality structure, the space of exterior forms is obtained
(quantum) antisymmetrizing the space of tensors and the quotient
calculus on the sphere is nontrivial ({\it i.e.} it is not
$0$--dimensional). The twisted limit of the multiparametric
calculus on the plane and of the multiparametric bicovariant
calculus on the corresponding orthogonal quantum group was first
studied in \cite{Aschieri}. In Section 3 we give a self contained
exposition of the Haar functional on twisted spheres and use it in
Section 4 to define an integral on forms. Since (differently from
the case of standard deformations) the Haar functional is a trace,
the integral on forms turns out to be a closed graded trace: we
thus obtain a {\it cycle} on the algebra. A cycle defines a cyclic
cocycle, a basic tool in noncommutative geometry \`a la Connes. In
Section 5 we define the pairing (metric) between tensorfields and
exploit the construction of forms as antisymmetrized tensors in
order to define the $*$--Hodge operator. Explicit expressions for
this operator and its properties are given. The $*$--Hodge
operator is shown to coincide with the one defined in
\cite{ConnesDuboisViolette}. Exterior forms, integration theory
and $*$--Hodge operator are fundamental ingredients for the
construction of field theories and gauge theories.

In the last Section we apply these tools to even
spheres. We make use of the projector introduced in
\cite{ConnesLandi} and \cite{ConnesDuboisViolette} which defines
the instanton bundle and we
compute its charge as the Chern--Connes pairing with the character
of the cycle defined in Section 4.

\bigskip
%%%%%%%%%%%%%%%%%%%%%% 1 %%%%%%%%%%%%%%%%%%
\thispagestyle{empty}
\section{Orthogonal Multiparametric Quantum Groups, Planes and Spheres}
The orthogonal multiparametric quantum groups are freely generated
by the $N^2$ matrix elements $\T{a}{b}$ (fundamental
representation) and the identity $1$, modulo the quadratic $RTT$
relations and the orthogonality relations discussed below.  The
noncommutativity is controlled by the multiparametric $R$-matrix
$R_{q,r}$
\begin{equation}
\Rmia{ab}{ef} \T{e}{c} \T{f}{d} = \T{b}{f} \T{a}{e} \Rmia{ef}{cd}
\label{RTTsosp}
\end{equation}
which satisfies the quantum Yang-Baxter equation. The
multiparametric $R_{q,r}$ matrix is obtained from the
uniparametric one $R_r$, via the transformation
\cite{Schirrmacher} (we follow the notations of \cite{Aschieri}):
$
R_{q,r}=F^{-1}R_rF^{-1}\,
$
where $(F^{-1})^{ab}_{~~cd}$ is a diagonal matrix
in the index couples $ab$, $cd$
\begin{equation}
F^{-1}\equiv diag (\sqrt{{r \over q_{11}}} ,
\sqrt{{r \over q_{12}}} , ... ~ \sqrt{{r \over q_{NN}}})
\label{effe}~,
\end{equation}
and where the complex parameters $q_{ab}$, $a,b=1,...N$ satify the
following relations
\begin{equation}
q_{aa}=r~,~~q_{ba}={r^2 \over q_{ab}}~,~~
q_{ab}={r^2 \over q_{a\bp}}={r^2 \over q_{\ap b}}=q_{\ap\bp}~,
 \label{qab2}
\end{equation}
in  the last equality we defined primed indices as $\ap \equiv
N+1-a$. Relations (\ref{qab2}) also imply   $q_{a\ap}=r$,
therefore the $q_{ab}$ with $a < b \leq {N\over 2}$ give all the
$q$'s. One can also easily show that the non diagonal elements of
$R_{q,r}$ coincide with those of $R_r$. The matrix $F$ satisfies
$F_{12}F_{21}=1$ {\it i.e.}
$F^{ab}{}_{ef}F^{fe}{}_{dc}=\delta^a_c\delta^b_d $, the quantum
Yang-Baxter equation $F_{12}F_{13}F_{23} =F_{23}F_{13}F_{12}$ and
the relations $(R_r)_{12}F_{13}F_{23}=F_{23}F_{13}(R_r)_{12}$.
Notice that for $r=1$ the multiparametric $R$ matrix reduces to
$R=F^{-2}$. Let $\Rh$ be the matrix defined by $\Rhat{ab}{cd}
\equiv \Rmia{ba}{cd}$, then the multiparametric $\Rh_{q,r}$ is
obtained from $\Rh_r$ via the similarity transformation
\begin{equation}
\Rh_{q,r}=F\Rh_rF^{-1}~;
\end{equation}
the characteristic equation and the projector decomposition
of $\Rh_{q,r}$ are therefore the same as in the uniparametric case:
\begin{equation}
(\Rh-rI)(\Rh+r^{-1}I)(\Rh- r^{1-N} I)=0 \label{cubic}
\end{equation}
{\vskip -1cm}
\begin{equation}
\Rh=r P_S - r^{-1} P_A+ r^{1-N}P_0  \label{RprojBCD}
\end{equation}
{\vskip -.3cm}
with
\begin{equation}
\begin{array}{ll}
&P_S={1 \over {r+\rminus}} [\Rh+\rminus I-(\rminus+
r^{1-N})P_0]\\
&P_A={1 \over {r+\rminus}} [-\Rh+rI-(r- r^{1-N})P_0]\\
&P_0= Q_N(r) K\\
&Q_N(r) \equiv (g_{ab} g^{ab})^{-1}={{1-r^{-2}} \over
{(1- r^{-N})(1+r^{N-2})}}~,~~~~
K^{ab}_{~~cd}\equiv g^{ab} g_{cd}\\
&I=P_S+P_A+P_0
\end{array}
\label{projBCD}
\end{equation}
Orthogonality of $T$ reads
\begin{equation}
g^{bc} \T{a}{b}  \T{d}{c}= g^{ad} 1~,~~~
g_{ac} \T{a}{b}  \T{c}{d}=g_{bd} 1 \label{Torthogonality}
\end{equation}
\noi where $g^{ab}=g_{ab}=\delta_{ab'}$. The consistency of
(\ref{Torthogonality}) with the $RTT$ relations is due to
\begin{equation} g_{ab} \Rhat{bc}{de} = \Rhatinv{cf}{ad} g_{fe}
~~~,~~~~~ \Rhat{bc}{de} g^{ea}=g^{bf} \Rhatinv{ca}{fd}
~.\label{crc2}
\end{equation}
\noi These identities
hold also for $\Rh \rightarrow \Rh^{-1}$.
A multiparametric determinant
is defined by det${}_{q,r}T=\epsilon_{q,r}^{i_1\ldots i_N}
\T{1}{i_1}\ldots \T{N}{i_N}$
cf. $\!$(\ref{epsilonqr}), it satisfies
$($det${}_{q,r}T)^2=1$. Imposing also the relation
det${}_{q,r}T=1$ we obtain the special orthogonal quantum group
$SO_{q,r}(N)$.
\sk
The multiparametric orthogonal quantum plane
%
%$\Re_{q,r}^{N}$
%
is the algebra
freely generated by the
elements $x^a$ with commutation relations
\begin{equation}
\PA{ab}{cd} x^c x^d=0 \label{PRTT13}~.
\end{equation}
The element
$
c=x^ag_{ab}x^b
$
is central and
imposing the extra relation $c=1$ we obtain the
multiparametric orthogonal quantum sphere.
%
%$\,{\Sqr}^{\!\!\!\!\!\!\!N-1}\equiv\Re_{q,r}^{N}/[c-1]\;$ (here
%$[c-1]$ is the ideal generated by the relation $c=1$).
\sk
The costructures of the orthogonal multiparametric quantum
groups have the same form as in the uniparametric case:
the coproduct
$\D$, the counit $\epsi$ and the antipode $S$ are given by
$
\D(\T{a}{b})=\T{a}{b} \otimes \T{b}{c}
~,~~\epsi (\T{a}{b})=\delta^a_b~,~~
S(\T{a}{b})=g^{ac} \T{d}{c} g_{db}
$
the coaction on the quantum plane and sphere reads
$
\delta(x^a)=T^a_{~b}\otimes x^b   .
$
\sk
\noi{\bf{Differential Calculus}}\\
There are only two $SO_{q,r}(N)$-covariant first order differential
calculi on the quantum plane such that any $1$-form can be uniquely
written as sum of functions on the quantum plane times the basic
differentials $dx^i$: $f_i dx^i$ {\it i.e.} such that the bimodule
of $1$-forms is generated as a left (or right) module by the
differentials $dx^i$. The deformed commutation relations are
\cite{Carow-Watamura:1991zp}
(the multiparametric case appeared in \cite{Scarfone})
\begin{equation}
 x^adx^b=  r\Rha^{ab}_{~~cd} dx^c x^d~~,
\label{xdx}
\end{equation}
the other calculus is obtained by replacing $r\Rha^{ab}_{~~cd}$ with
$r^{-1}\Rha^{-1}{}^{ab}_{~~cd}$ (so that for $r=1$ the calculus is unique).
The exterior differential
$d$ by definition satisfies the Leibniz rule and therefore, from (\ref{xdx})
it follows that the algebra of exterior forms is generated by $x^a$ and $dx^a$
modulo the ideal generated by the relations (\ref{PRTT13}), (\ref{xdx})
and
\begin{equation}
\label{relations}
dx^a \, dx^b =-r {\hat R^{ab}_{~cd}}{}^{\,} dx^c{}^{\,} dx^d~~;
\end{equation}
recalling  (\ref{RprojBCD}) and (\ref{projBCD}) this relation is equivalent to
${P_S}^{ab}_{~~cd}dx^c {}^{\,}dx^d={P_0}^{ab}_{~~cd}dx^c{}^{\,}
dx^d=0$.
Partial derivatives can be defined so that $da = dx^c\partial_c(a)$.
They satisfy the deformed Leibniz rule and the commutation relations
\begin{equation}
\partial_c x^b= \delta^b_c 1 \,+
r \Rha^{be}_{~\;cd} x^d \,\partial_e~~~,~~~~~
 \PA{ab}{cd} \partial_b \partial_a=0~.
\end{equation}

\begin{remark}\label{forms_from_tensorproduct}{\rm
The space of $2$-forms defined by (\ref{relations}) is equivalent to the
space of $2$-forms defined using the wedge product
(here $\mu$ is an arbitrary coefficient)
\begin{eqnarray}
dx^a\wedge dx^b&\equiv& \mu P_A^{ab}{}_{~~cd}{}^{\,}dx^c\otimes dx^d
\label{dwd}\\
&=&dx^a\otimes dx^b-(I-\mu P_A)^{ab}_{~~cd}dx^c\otimes dx^d
%\equiv dx^a \wedge dx^b-\Lambda^{ab}_{~~cd}dx^c\otimes dx^d
\end{eqnarray}
Notice however that forall $\mu$,
\begin{equation}
\label{Lambda}
\Lambda\equiv (I-\mu P_A)={\mu\over r+r^{-1}}[-\Rh+\mu^{-1}(\mu r -r- r^{-1} )
I-(r-r^{1-N})P_0]
\end{equation}
does not satisfy the braid equation.
This situation differs from that of the $GL_{q,r}(N)$-covariant
plane \cite{Wess:1991vh} and from that of
exterior forms on quantum groups as described in \cite{Woronowicz}.
There the corresponding $\Lambda$ matrix satisfies the braid equation
so that the space of $k$-forms (not only that of $2$-forms) can be
defined as the space of quantum
antisymmetric tensorfields.
We have
\begin{equation}
\label{W}
dx^{i_1}\wedge\ldots dx^{i_k}=W(dx^{i_1}\otimes\ldots
dx^{i_k})=W^{i_1\ldots i_k}_{~~j_1\ldots j_k} dx^{j_1}\otimes
dx^{j_k}\;,
\end{equation}
where the numerical coefficients $W_{1\ldots k}= W^{i_1\ldots
i_k}_{~~j_1\ldots j_k}$ give the alternating sum of $k!$ addends,
these addends corresponding to the $k!$ permutations of $k$
elements. Since $\Lambda$ is not in general a representation of
the permutation group each permutation must be expressed via a
minimal set of nearest neighbour transpositions, each
transposition is then represented via $\Lambda$. A recursion
relation for $W$ is
\begin{equation}
\label{wedge2}
W_{1\ldots k} = {\cal I}_{1\ldots k} W_{1\ldots
k-1} ,
\end{equation}
{\vskip -.4cm} \noi where ${\cal I}_{1\ldots k} = I -
\Lambda_{k-1,k} + \Lambda_{k-2,k-1} \Lambda_{k-1,k}   \ldots
-(-1)^k \Lambda_{12} \Lambda_{23} \cdots \Lambda_{k-1,k}$ and
$W^i{}_j = {\cal I}^i{}_j = \delta^i_j$. A recursion relation for
${\cal I}$ is
\begin{equation}
{\cal I}_{1\ldots k} = I - {\cal I}_{1\ldots k-1} \Lambda_{k-1,k}
\label{wedge4}~~.
\end{equation}

A different space of $2$-forms on the quantum orthogonal plane
is defined imposing only the relation
${P_S}^{ab}_{~~cd}dx^c {}^{\,}dx^d=0^{\,}$, cf. \cite{FRT};
this relation (and only this one for generic $r$),
is implied by the  wedge product
%[recall (\ref{RprojBCD}),(\ref{projBCD})]
$dx^a\wedge dx^b\equiv dx^a\otimes dx^b -
r^{-1}\Rh^{ab}_{~~cd}dx^c\otimes dx^d$. In this case
$\Lambda=r^{-1}\Rh$ satisfies the braid equation and the space of
exterior forms can be constructed as in (\ref{W})--(\ref{wedge4}).
A differential calculus on the quantum  orthogonal plane with this
exterior algebra is studied in \cite{Scarfone}; it is obtained as
a bicovariant calculus (\cite{Woronowicz}) on the inhomogeneous
orthogonal quantum group $ISO_{q,r}(N).\!\!\!$ \fidi}\end{remark}
%\medskip

\begin{remark}
\rm{In the next sections we consider the $r\rightarrow 1$ limit of the
exterior algebra defined by (\ref{relations}), equivalently (\ref{dwd}).
In this limit the
$\Lambda$ matrix (\ref{Lambda}) with  $\mu=2$ equals the
$\Rh_{q,r=1}$ matrix  (cf. last paragraph in Remark 1), it satisfies
the braid equation and it also squares to the identity matrix,
therefore defining a representation of the permutation group.
We can thus construct the space of exerior forms
as in (\ref{W})--(\ref{wedge4}). Here we show that the exterior algebra
Im($W$) obtained as the image of the antisymmetrizer map $W$
is isomorphic to the exterior algebra $\Omega$ freely generated by the
elements $dx^i$ modulo the ideal generated by the relations
(\ref{relations}).
Indeed consider the surjection $\phi\,: ~\Omega\longrightarrow {\rm{Im}}W$
defined by $\phi([{\cal T}])={1\over k!}W({\cal T})$ where
${\cal T}$ is a polynomial of degree $k$ in the $dx^i$, and $[\cal T]$
denotes the equivalence class under (\ref{relations}).
$\phi$ is well defined, it is also injective because
$[{\cal T}]= [{1\over k!}W({\cal T})]$;
this last relation is proven observing that
$\Lambda^{ab}_{~~cd}[dx^c dx^d]= -[dx^a dx^b]$. A similar
argument shows that $\Omega\cong$ Im($W$) also for the
$GL_{q,r\not=1}(N)$--convariant plane.
\fidi}\end{remark}

Generalizing \cite{Fiore:1994vk} to the multiparametric case,
because of (\ref{relations}) and of the specific expression
of $\Rh_{q,r}$, we have that any monomial
$dx^{i_1}dx^{i_2}\ldots dx^{i_p}$ can be rewritten
as sum of monomials  $dx^{j_1}dx^{j_2}\ldots dx^{j_p}$ with
$j_1<j_2<\ldots<j_p$ so that the graded differential algebra
of exterior forms on the quantum plane
%
%$\Omega(\Re_{q,r}^N)=\bigoplus_{k=0}^N \Omega_k$
%(with $\Omega_0=\Re_{q,r}^N$)
%
has dimension $2^N$ as in the classical case.
In particular every $N$-form is proportional to the volume
form\footnote{We show in Subsection 5.2 that $V_N$ is real.}
\begin{equation}
\label{VN}{}~~~~~~~~~~~V_N=i^{[{N\over 2}]}\,dx^1{}^{}dx^2\ldots dx^N~~,~~~~~~~~~
\sma{$[{N\over2}]\equiv\,$integer part of ${N\over 2}$}~.
\end{equation}
The epsilon tensor is defined by
\begin{equation}
\epsilon_{q,r}^{i_1\ldots i_N} dx^1dx^{2}\ldots dx^N
= dx^{i_1}dx^{i_2}\ldots dx^{i_N}
\label{epsilonqr}~.
\end{equation}
One can show that $x^a V_N= r^N V_N x^a\,$.\footnote{{\sl Proof.}
In order to show $\epsilon_{q,r}^{i_1\ldots i_N} \Rh^{a 1}_{~i_1
a_1}\Rh^{a_1 2}_{~\,i_2 a_2}\ldots \Rh^{a_{N-1} N}_{~\,i_N
a_N}=\pm\delta^a_{~a_N}$ apply the $SO_{q,r}(N)$-coaction to
$V_N$, recall that $($det$T)^2=1$ and consider the $N\times N$
representation of $\T{c}{d}$ given by $\Rh^{ac}_{~\,db}$. Finally
the plus sign is singled out going to the commutative limit
$r=q_{ab}=1$.} If we extend the quantum plane algebra including
the generator $c^{-1}$, the exterior differential is then given by
the $1$-form \cite{Steinacker}
$
\omega={r^2\over{r+1}}c^{-1}dc
%=r{c^{-{1\over 2}}}d c^{1\over 2}
%=dc^{1\over 2}\; c^{{-{1\over 2}}}
$
as follows
$
d\th={1\over{1-r}}\,[\omega,\th]_{\pm}
$
where we use the commutator if $\theta$ is an even form, the
anticommutator if $\theta$ is odd. Notice that $d$ is an inner
differential only if $r\not =1$. The drop of dimension discussed
in the introduction is related to this property of the exterior
differential. While this aspect may seem a trivialization (from an
outer $d$ with $r=1$ to an inner $d$ with $r\not=1$) it also hints
that the geometry is highly noncommutative, indeed $d$ and the
partial derivatives are finite difference operators for $r\not
=1$.

\sk It is natural to study how the calculus on the $N+1$
dimensional quantum plane induces a calculus on the $N$
dimensional sphere. As in the commutative case we define the
exterior algebra on the sphere as the quotient of the exterior
algebra on the plane modulo the differential ideal generated by
the relation $c=1$. Since $c$ is not central in the differential
algebra, {\it i.e.} $c\ dx^a=r^2 dx^a c$ and $x^a dc=dc\ x^a
+(1-r^{-2})c\ dx^a$, we immediately have that $dx^a=0$ if
$r\not=1$. We conclude that in the $r\not=1$ case the quotient
calculus on the sphere is trivial.

\sk
\noi{\bf Real Forms}\\ All real forms of (uniparametric)
orthogonal quantum groups and their quantum spaces are studied in
\cite{Twietmeyer:1992mj}, (see also \cite{Dobrev}). Here we focus
on the compact form $SO_{q,r}(N,\Re)$ and on the multiparametric
Euclidean quantum plane $\Re_{q,r}^{N}$ and sphere
${\Sqr}^{\!\!\!\!\!\!\!N}$. These are given by the conjugation
\begin{equation}
(\T{a}{b})^{\*}=g_{ea}\T{e}{f}g^{bf}~~,~~~~  (x^a)^{\*}=g_{ea}x^e~~,
\label{conjugation}
\end{equation}
that is compatible with the quantum group, plane and sphere
defining relations and with the coaction
$
\delta(x^a)=T^a_{~b}\otimes x^b
$
if ${\bar{R}}^{ab}_{~cd}=
R^{dc}_{~ba}$, {\it {i.e.}}
$q_{ab} {\bar q}_{ab}=r^2, r \in{\Re}\,.$
Conjugation (\ref{conjugation}) however, for $r\not=1$, is not compatible
with the differential calculus on the
quantum plane in the sense that (\ref{xdx}) implies
$(d x^a)^\*\not= d (x^{a\,\*})$.
Also, the conjugated partial derivatives
$\pa^\*_a$ are not linear combinations of the $\pa_a$'s. Rather
$(d x^a)^{\*}$ and $\pa^\*$ generate the other calculus on
the quantum orthogonal plane [cf. line after (\ref{xdx})].
In \cite{Ogievetsky:1992qp}, $(d x^a)^\*$ and $\pa^\*_a$
are expressed nonlinearly in terms of the $x,dx,\pa$ algebra.
We just mention that on the other hand the conjugations that give the
signatures $n,m$ with $n+m=N$, and $n-m=0,1,2$ give a real
differential calculus $(dx^a)^\*=d(x^{a\,\*})$ \cite{Scarfone}.
\sk
\noi{\bf{Integration}}\\
{}Generalizing \cite{Steinacker}
to the multiparametric case\footnote{Hint: use induction on the number
of deformation parameters. The positivity property
$h([x^{i_1}x^{i_2}\ldots x^{i_{2n+1}}])=
\overline{\,h([{x^{i_{2n+1}}}^\*
\ldots {x^{i_2}}^\*{x^{i_1}}^\*])\,}$
holds  because of $g_{ij}\in \Re$, (\ref{conjugation}),
${\overline{{P_A}^{ef}_{~ih}}}= {P_A}^{ih}_{~ef}~\rm{and}$
${P_A}^{ih}_{~ef}g_{ji}g_{kh}= {P_A}^{\ell m}_{~kj}g_{me}g_{\ell f}$,
cf. (\ref{Torthogonality}). Contrary to \cite{Steinacker} reality of
$\hat R$ is not needed.}
, we
obtain that there exists a unique (normalized) integral
of functions on the multiparametric sphere $\Sqr^N$ such that
it is invariant under the $SO_{q,r}(N,\Re)$ coaction and
it is analytic in $r-1$ and $q_{ab}-r$.
We use the notation $h(f)$ for the integral
(Haar functional) of $f\in \Sqr^N$.
On the elements
$[x^{i_1}x^{i_2}\ldots x^{i_{p}}]\in \Sqr^N$
(the square brakets denote the equivalence class
w.r.t. $\!$the relation $c=1$) we have
$h([x^{i_1}x^{i_2}\ldots x^{i_{2n+1}}])=0$ and
$h([x^{i_1}x^{i_2}\ldots x^{i_{2n}}])=\lambda_n
\Delta^n(x^{i_1}x^{i_2}\ldots x^{i_{2n}})$ with
$\Delta=g^{ji}\pa_i\pa_j$ and $\lambda_n$ a
proportionality factor depending only on $n$ and $r$.
The Haar functional on  $\Sqr^N$ has the
following reality, positivity and
quantum cyclicity  properties:
$$\overline{h(f)}=h(f^\*)~~,~~~~h(f^\* f)\geq 0~~,~~~~
h(f g)=h(g \,D\!f)
$$
where $f,g,D{\!}f\in \Sq^N$.
The  map
$D$ is defined on the basic monomials $[x^{i_1}x^{i_2}\ldots x^{i_p}]$ as
$D[x^{i_1}x^{i_2}\ldots x^{i_p}]\equiv
[D^{i_1}_{~j_1}x^{j_1}_{\,}D^{i_2}_{~j_2}x^{j_2}\ldots
D^{i_p}_{~j_p}x^{j_p}]$, where $D^a{}_e\equiv g^{as}g_{es}$.
This map $D$ is then extended by linearity to all of $\Sqr^N$.
(It is easy to see that $D$ is well defined, in fact $D[c-1]=0$).
In the twisted limit $r\rightarrow 1$, we have
$D^a{}_e\rightarrow \delta^a{}_e$  and we obtain the  cyclic
property $h(f g)=h(g f)$. In the next section we give a self-contained
exposition of the Haar functional on twisted spheres, we also give
the explicit expression of the $\lambda_n$ coefficients.
\bigskip
\section{Twisted Spheres}
The twisted quantum Euclidean planes and spheres  $\Re_q^N$,
$\Sq^N$ are obtained considering the limit $r\rightarrow 1$ of the
corresponding multiparametric structures. We have ($a'=N+1-a$)
\begin{equation}
\Rh_q^{~ab}{}_{\,cd}=
q_{ab}\,\delta^a_{\;d}\delta^b_{\;c}~~~,~~~~~
g_{ab}=g^{ab}=\delta_{ab'}
\end{equation}
\begin{equation}
|q_{ab}|=1~,~~q_{aa}=q_{aa'}=1~,~~ q_{ab}=q_{a'b'} ~,~~
q_{ab}=q_{ba}^{-1}=q_{ab'}^{-1}\;.
\end{equation}
{}For each $p$ we have that
\begin{equation}
\label{prodotto_dei_q} \Pi_{i=1}^N q_{ip}= 1   \;.
\end{equation}
Explicitly, the twisted quantum Euclidean plane
% $\Re_q^N$
$q$-commutation relations are
\begin{equation}
\label{relations3}
x^ax^b = q_{ab} x^b x^a
\end{equation}
The twisted quantum Euclidean sphere is the quotient algebra
$\Sq^N=\Re_q^{N+1}/{\rm{I}}$, where $\rm{I}$ is the ideal generated by
the relation $c=1$.
\medskip
\begin{remark}\label{the_explicit_isomorphism}{\rm
If we specialize to $\Sq^4$ the only independent deformation
parameter is $q_{12}=q$. The explicit relations are
$2 (x^1x^5+x^2x^4) + (x^3)^2 = 1$, $[x^3,~.~]=0$ and
\begin{eqnarray*}
x^1 x^2 = q x^2 x^1 ~,~~ & x^1 x^4 = q^{-1} x^4 x^1  ~,~~& x^1 x^5
= x^5 x^1 ~,  \cr x^2 x^5 = q x^5 x^2  ~,~~ & x^4 x^5 = q^{-1} x^5
x^4 ~,~~& x^2 x^4 = x^4 x^2  ~.
\end{eqnarray*}
The explicit isomorphism with Connes-Landi sphere is given by
$\lambda=q$, $x^3=2t-1$, $x^1=\sqrt{2}\alpha$,
$x^5=\sqrt{2}\alpha^\*$, $x^2=\sqrt{2}\beta$,
$x^4=\sqrt{2}\beta^\*$.
More generally, setting $q_{ab}=e^{i\theta_{ab}}$ forall $a<a',b<b'$,
the twisted planes and shperes coincide with the $\theta$-deformed
ones defined in \cite{ConnesDuboisViolette}.
\fidi}
\end{remark}
\medskip
{\bf Integration of functions on the sphere}

\noi Let
$\Delta=g^{ji}\partial_i\partial_j$ be the Laplacian in
$\Re_q^N$, where the partial derivatives now satisfy, forall
$f\in\Re_q^N$,
\begin{equation}
\partial_s(x^af)=\delta^a_s f + q_{as}
x^a\partial_s f~~,~~~~
\partial_a \partial_b=q_{ab}\partial_b \partial_a~.
\label{relations4}
\end{equation}
A straightforward computation shows that (no sum on $k$)
\begin{equation}
\label{laplaciano_intermedio} \Delta ~ x^kx^{k'} =
 x^kx^{k'}~\!\Delta\, + 2 + 2 x^k\partial_k +2 x^{k'}\partial_{k'}\;,
\end{equation}
so that
\begin{equation}
\label{laplaciano_metrica} \Delta ~ c = c~\Delta + 2N + 4
x^j\partial_j \;.
\end{equation}
Before introducing
the Haar functional we have to show the following lemma.
\begin{lemma} For each $n>0$ we have that
\begin{equation}
\label{Delta_Metrica} \Delta^{n+1}(c x^{i_1}\ldots x^{i_{2n}})=
2(n+1)(N+2n) \Delta^n(x^{i_1}\ldots x^{i_{2n}})\;.
\end{equation}
\end{lemma}
{\it Proof}. We first show that forall $n$ and for
each $0\leq k\leq n$
we have
\begin{equation}
\label{recurrencemm}
\begin{array}{l}
 a_k \Delta^n(x^{i_1}\ldots x^{i_{2n}}) +
\Delta^{n-k+1}c\Delta^k(x^{i_1}\ldots x^{i_{2n}})=~~~~~~~~~~~~~~~~~~~~~~~~~~~~~~~~~~~~~\\
~~~~~~~~~~~~~~~~~~~~~~~~~~~~~~~~~~~ a_{k+1} \Delta^n(x^{i_1}\ldots x^{i_{2n}}) +
\Delta^{n-k}c\Delta^{k+1}(x^{i_1}\ldots x^{i_{2n}})~~~
\end{array}
\end{equation}
with $a_k= 2k(N+4n)-4k(k-1)$.
Indeed recalling (\ref{laplaciano_metrica})  and observing that
$x^j\partial_j(x^{i_1}\ldots x^{i_{2(n-k)}}) = 2(n-k)
x^{i_1}\ldots x^{i_{2(n-k)}}\;
$ the l.h.s. of (\ref{recurrencemm})
equals
$$
(a_k + 2 N
+8(n-k))\Delta^n(cx^{i_1}\ldots x^{i_{2n}}) +
\Delta^{n-k}c\Delta^{k+1}(x^{i_1}\ldots x^{i_{2n}})\;.
$$
Formula (\ref{recurrencemm}) is then shown by verifying that
$a_{k+1}=a_k+2N+8(n-k)$. Property (\ref{Delta_Metrica}) then
follows by choosing in (\ref{recurrencemm}) $k=0$ and $k=n$ and observing
that $\Delta^{n+1}(x^{i_1}\ldots
x^{i_{2n}})=0$. \fidi
\medskip
Recalling that $(x,a)_n=x(x+a)\ldots (x+(n-1)a)$ and $(x,a)_0=1$,
let us define $h_\Re:\Re_q^{N+1}\rightarrow\Co$ as the linear map
that on monomials is given by
\begin{equation}
\label{functional_on_R} h_\Re (x^{i_1}\ldots x^{i_{2n}}) =
\lambda_n\Delta^n(x^{i_1}\ldots x^{i_{2n}})~~,~~~~
\lambda_n=\frac{1}{2^n \,n!\,(N,2)_n}\;.
\end{equation}
\bigskip
\begin{proposition}
\label{Haar} Let $f\in \Re_q^{N+1}$ and let $[f]\in \Sq^N$ be its
equivalence class. The linear functional $h([f]) \equiv h_\Re(f)$
is well defined on $\Sq^N$ and satisfies the following properties
(we omit to denote the equivalence class) :
\begin{itemize}
\item[{\rm a)}] $h(1)=1$.
\item[{\rm b)}] $h(f g)=h(g f)$, $~~~~~~~~~~~~~~~~~~~~$for each $f,g\in\Sq^N$.
\item[{\rm c)}] $1\otimes h=(\id\otimes h)\sma{$^{\,}\circ^{\,}$}\delta$
~~~~~~~~~~~~~~{{(where here, $1$ is the identity in $SO_{q,r}(N)\,$)}}~.
\item[{\rm d)}] $\overline{h(f)}=h(f^\*)$ and $h(f^\*f)>0$, for each
$f\in\Sq^N\,$ (reality and positivity of $h$).
\end{itemize}
\end{proposition}
{\it Proof}.  From the definition of $\lambda_n$ and from Lemma
(\ref{Delta_Metrica}) we directly check that
\begin{eqnarray*}
h_\Re(x^{i_1}\ldots x^{i_{2n}}(c-1)) &=&
\lambda_{n+1}\Delta^{n+1}(x^{i_1}\ldots
x^{i_{2n}}c)-\lambda_n\Delta^n(x^{i_1}\ldots x^{i_{2n}})\cr
&=&(2(n+1)(N+2n)\lambda_{n+1}-\lambda_n)\Delta^n(x^{i_1}\ldots
x^{i_{2n}}) = 0\;,
\end{eqnarray*}
{\it i.e.} $h_\Re((c-1)\Re_q^{N+1})=0$ and $h$ is well defined.
Point a) is trivial. To prove point b) let us remark that if
$h(x^{i_1}\ldots x^{i_{2n}})\not = 0$ then for every index
$\ell\in \{i_1,\ldots,i_{2n}\}$
there is a companion index $\ell'\in \{i_1,\ldots,i_{2n}\}$, {\it i.e.}
$\{i_1,\ldots,i_{2n}\}=\{j_{1},j_{1}'\ldots j_n,j_n'\}$. It follows
\begin{eqnarray*}
h(x^{i_1}\ldots x^{i_{2n}}) = \Pi_{k=1}^{2n} q_{i_ki_{2n}}
h(x^{i_{2n}}x^{i_1}\ldots x^{i_{2n-1}}) &=& \Pi_{k=1}^{n}
q_{j_ki_{2n}}q_{j_k'i_{2n}} h(x^{i_{2n}}x^{i_1}\ldots
x^{i_{2n-1}})\cr &=& h(x^{i_{2n}}x^{i_1}\ldots x^{i_{2n-1}})\;.
\end{eqnarray*}
Since the calculus is covariant, the partial derivatives satisfy the following
property: $$ \delta\sma{$^{\,}\circ^{\,}$}\partial_\ell=
(S^{-1}(\T{k}{\ell})\otimes\partial_k)\sma{$^{\,}\circ^{\,}$}\delta~~.$$ The
Laplacian $\Delta$ is then $SO_q(N+1)$--invariant, {\it i.e.}
$(\id\otimes\Delta)\sma{$^{\,}\circ^{\,}$}\delta=\delta\sma{$^{\,}\circ^{\,}$}\Delta$, and point c) is
then proved.

In order to prove reality of $h$ we first
observe that $\Re_q^N$ can be linearly generated by the ordered
monomials
$x^{e_1}\ldots x^{e_k} x^{a_1}x^{a_1'}\ldots x^{a_s}x^{a_s'}$,
with $e_1\leq e_{2}\ldots \leq e_k$ and where the $e_i$ indices
do not have a companion $e_i'$ index.
The action of  $x^k\partial_k$ (no sum on $k$)
on these monomials doesn't depend on the $q_{ab}$, so that, thanks to
(\ref{laplaciano_intermedio}), the value of $h$ on these monomials
equals the classical value. Reality then easily follows from the
reality of the classical integral. Analogously, to prove
positivity it is enough to verify that if $f$ is an
ordered polynomial the terms of $f^\*f$ which contribute to $h$
are automatically ordered monomials. \fidi
\medskip
We call $h$ the {\it Haar functional} on $S^N_q$.
We have seen that $h$
on the ordered monomials
$x^{e_1}\ldots x^{e_k} x^{a_1}x^{a_1'}\ldots x^{a_s}x^{a_s'}$
equals the commutative integral.

\bigskip
\section{Calculus on $\Re_q^N$ and $\Sq^N$}
The graded differential algebra  $\Omega(\Re_q^N)=\bigoplus_{k=0}^N
\Omega_k(\Re_q^N)$ , with $\Omega_0(\Re_q^N)=\Re_q^N$,
is the $r\rightarrow 1$ limit
of the multiparametric one. As shown in Remark 2, for $r=1$
we can consider  $\Omega(\Re_q^N)$ as the space  of completely
$q$-antisymmetrized tensors. The explicit relations are
(\ref{relations3}) and
\begin{equation}
\label{relations2}  dx^a x^b= q_{ab}  x^b dx^a ~~,~~~~ dx^a\wedge
dx^b= -q_{ab_{\,}} dx^b\wedge dx^a ~.
\end{equation}
Conjugation (\ref{conjugation}) is now compatible with the
differential calculus, we have $x^{a\,\*}=x^{a'}$ and $(d x^a)^\*=
d (x^{a\,\*})$. The volume form (\ref{VN}) is now central [this
result follows also directly from (\ref{prodotto_dei_q})].

\medskip
We now define a differential graded algebra on $\Sq^N$. Let $J=J_1
+ J_2\subset\Omega(\Re_q^{N+1})$, where
$$
J_1 = \{(c-1)\omega,\ \omega\in\Omega(\Re_q^{N+1})\}~~,~~~~ J_2 =
\{\omega\in\Omega(\Re_q^{N+1}), \  \omega\wedge dc=0\} \;.
$$
Because $c$ and $dc$ are central, both $J_1$ and $J_2$ are graded
ideals. Moreover it is easy to verify that $d(J)\subset J$ and
$\delta(J)\subset SO_q(N+1)\otimes J$, {\it i.e.} $J$ is a
differential ideal and left coideal. We have that
$$
\Omega(\Sq^N) \equiv \Omega(\Re_q^{N+1})/J\equiv \bigoplus_{k=0}^N
\Omega_k(\Sq^N)
$$
is a left covariant differential graded algebra with
$\Omega_0(\Sq^N)=\Sq^N$. In the following we denote with
$[\omega]\in\Omega_k(\Sq^N)$ the equivalence class of
$\omega\in\Omega_k(\Re_q^{N+1})$; we have that $d[\omega]=[d\omega]$.

Let $\omega_k\in\Omega_N(\Re_q^{N+1})$ be defined by
\begin{equation}\label{omega}
\omega_k = \frac{1}{N!}
i^{[\frac{N+1}{2}]}\epsilon_{{}_{^{\!\!(q^{\!-\!1\!})}}}{}_{s_1\ldots
s_Nk}dx^{s_1}\wedge\ldots dx^{s_N}\;,
\end{equation}
where  ${{\epsilon_{{}_{^{\!\!(q^{\!-\!1\!})}}}}}$
 is the epsilon tensor ${{\epsilon_{{}_{^{\!\!(q^{\!-\!1\!}),r=1}}}}}$
[c.f. (\ref{epsilonqr}) and (\ref{Wem1})].  Thanks to
(\ref{eem1bis=W}) we have that $\omega_k \wedge dx^\ell =
\delta^\ell_k V_{N+1}$ and therefore
$$
 x^k \omega_k \wedge
dc=x^k\omega_k\wedge(x^a g_{ab}dx^{b}+ dx^a
g_{ab}x^{b} )
= 2x^k\omega_k \wedge dx^{a}g_{ab} x^{b} =
2 x^k V_{\!N+1} g_{kb}x^{b} = 2c V_{\!N+1}.
$$
{}From this formula, using that
on a commutative sphere of unit radius we have $dc/2=d\,(g_{ab}x^a x^b)^{1\over 2}$, we read off the volume form on
$\Omega(\Sq^N)$:
$$
\V_N=[x^k\omega_k] \;.
$$
Any $N$ forms on the sphere can be expressed in terms of $\V_N$
as follows.
Let $\omega\in\Omega_N(\Re_q^{N+1})$ so that
 $\omega\wedge dc/2$ is proportional to  $V_{N+1}$, we set
$\omega\wedge
dc/2=f_\omega V_{N+1}$. We then have
$$
c\omega \wedge dc = 2 f_\omega c V_{N+1}= f_\omega x^k
\omega_k\wedge dc~~~,
$$
and therefore $c\omega-f_\omega x^k\omega_k \in J_2$. Since
$\omega-c\omega\in J_1$, we obtain that
$[\omega]=[f_\omega x^k\omega_k]$ $=[f_\omega]\V_N$.

\medskip
We are now ready to define an integral on $N$-forms:
\begin{equation}
\label{integral} \int [\omega] = \int [f_\omega]\V_N =
h([f_\omega])\;,
\end{equation}
where $\omega\wedge dc/2=f_\omega V_{N+1}$. This integral verifies
the following Stokes' Theorem:
\begin{proposition}\label{stokes}
$$
\int d[\theta] = 0 ~~~~~\forall ~ [\theta]\in\Omega_{N-1}(\Sq^N)~.
$$
\end{proposition}
{\it Proof}. Let $\theta=dx^{i_1}\wedge\ldots
dx^{{i_{N-1}}}f_{i_1\ldots
i_{N-1}}\in\Omega_{N-1}(\Re_q^{N+1})$.
By direct computation we obtain that
\begin{eqnarray*}
d\theta \wedge dc &=&
(-1)^{N-1}2 dx^{i_1}\wedge\ldots
dx^{{i_{N+1}}}x^\ell\partial_{i_N}
f_{i_1\ldots
i_{N-1}}g_{i_{N+1}\ell} \cr
&=& (-1)^{N-1} 2 \epsilon_q^{i_1\ldots i_{N+1}}
x^\ell\partial_{i_N} f_{i_1\ldots
i_{N-1}}g_{i_{N+1}\ell}V_{N+1}\;.
\end{eqnarray*}
By using the definition (\ref{integral})
of the integral on $N$--forms we have that if $f_{i_1\ldots i_{N-1}}$
 is odd in the coordinates $x^i$ then $\int d[\theta] =0$, while if
$f_{i_1\ldots i_{N-1}}$ is even and of degree $2n$ in the coordinates
$x^{i}$ we have
\begin{eqnarray*}
\int d[\theta] &=& (-1)^{N-1}\epsilon_q^{i_1\ldots
i_{N+1}} h(x^\ell\partial_{i_N} f_{i_1\ldots i_{N-1}})g_{\ell
i_{N+1}}\cr
&=&(-1)^{N-1}\epsilon_{q_{{}_{_{}}}}^{i_1\ldots
i_{N+1}} \lambda_n\Delta^n(x^\ell\partial_{i_N} f_{i_1\ldots
i_{N-1}})g_{\ell i_{N+1}}\cr
&=&(-1)^{N-1}\epsilon_{q_{{}_{_{}}}}^{i^{^{}}_1\ldots
i_{N+1}} \lambda_n\Delta^{n-1}(2\partial_t g^{t\ell}+
x^\ell\Delta) \partial_{i_N} f_{i_1\ldots i_{N-1}}g_{\ell
i_{N+1}}\cr
&=&(-1)^{N-1}\epsilon_{q_{{}_{_{}}}}^{i^{^{}}_1\ldots
i_{N+1}}
\lambda_n\Delta^{n-1}(2\partial_{i_{N+1}}\partial_{i_N}f_{i_1\ldots
i_{N-1}}+ g_{i_{N+1}\ell}x^\ell\Delta
\partial_{i_N} f_{i_1\ldots i_{N-1}})\cr
&=&(-1)^{N-1} \epsilon_{q_{{}_{_{}}}}^{i^{^{}}_1\ldots
i_{N+1}} \lambda_n\Delta^{n-1} g_{i_{N+1}\ell}x^\ell\Delta
\partial_{i_N} f_{i_1\ldots i_{N-1}}\;,
\end{eqnarray*}
where we used the relation $[\Delta
,x^\ell]=2 g^{t\ell}\partial_t$ and the $q$--antisymmetry of
$\epsilon_q$--tensor, {\it i.e.} $\epsilon_q^{i_1\ldots
i_Ni_{N+1}}=-q_{i_Ni_{N+1}} \epsilon_q^{i_1\ldots i_{N+1}i_N}$
together with
$\partial_{i_{N+1}}\partial_{i_{N}}=
q_{i_{N+1},i_{N}}\partial_{i_{N}}\partial_{i_{N+1}}$. Since
$\Delta\partial_k=\partial_k\Delta$ we can repeat the above
argument and obtain
$$
\int d[\theta] = (-1)^{N-1} \epsilon_q^{i_1\ldots
i_{N+1}} \lambda_ng_{i_{N+1}\ell}x^\ell\Delta^n
\partial_{i_N} f_{i_1\ldots i_{N-1}} = 0\;,
$$
because $\deg (\partial_{i_N}f_{i_1\ldots i_N})=2n-1$. \fidi
\medskip
The integral (\ref{integral}) has also the property, for each $[a]\in\Sq^N$ and $[\omega]\in\Omega_N(\Sq^N)$,
\begin{eqnarray*}
\int [a\omega]&=&\int [a f_\omega]\V_N =
h([af_\omega])=h([f_\omega a])\cr
&=&\int [f_\omega a]\V_N = \int
[f_\omega]\V_N [a] = \int [\omega a]~~.
\end{eqnarray*}
Following the proof of Proposition III.$4$ ($1\implies 3$) of
\cite{ConnesBook} we can conclude that the integral $\int$ is a
closed graded trace on $\Omega(\Sq^N)$. We summarize the results
of this section in
the following proposition (for a definition of cycle see Section 6):
\medskip
\begin{theorem}
\label{cycle} {\rm ($\Omega(\Sq^N),d,\int$)} is a cycle.$~\;$\fidi
\end{theorem}
\medskip

{}For future reference, we denote with $\tau$ the character of the cycle
($\Omega(\Sq^N),d,\int$),
\begin{equation}
\label{carattere} \tau(a_0,a_1\ldots a_N)=\frac{2^{[N/2]+1}[N/2]!}{i^{[N/2]}N!}
 \int a_0 da_1\ldots
da_N~~~~~~~,~~~ a_i\in\Sq^N\;.
\end{equation}
The normalization in this formula is chosen in order to fix the charge of the Bott
projector on the classical even sphere equal to $1$. Indeed if $p_{\!B}$ is
the Bott projector for $S^{2n}$ it can be shown that its Chern character
is $ch(p_{\!B})=1+\frac{i^n(2n)!}{2^{n+1}n!}\V_{2n}$ (see
\cite{VarillyBook}). Since the definition of the Haar measure in
Proposition \ref{Haar} and in (\ref{functional_on_R}) doesn't
contain $q$--factors, it is natural not to $q$--deform the
normalization of the character $\tau$.

%%%%%%%%%%%%%%%%%%%%%%%%%%%%%%%%%%%%%%%%%%%%%%%%%%%
\bigskip
\section{Hodge Theory}
\subsection{Hodge Theory on $\Re_q^N$}
We already observed that in the $r\rightarrow 1$ limit
the space of exterior forms is the image of the $q$-antisymmetrizer
$W$ introduced in (\ref{W}).
Moreover since in this case $\Lambda^{ab}_{~cd}\propto
{\Lambda_{{}^{_{\!(q{=1})}}}}^{\!\!\!ab}_{cd}$, every permutation
$W^{i_1\ldots i_k}_{~~j_1\ldots j_k}$ differs from the $q_{ab}=1$
permutation ${W_{{}^{_{\!(q{=1})}}}}^{\!\!\!i_1\ldots
i_k}_{\,j_1\ldots j_k}$ by at most a proportionality factor given
by a monomial in the $q_{ab}$'s. In particular (no sum on $i$'s)
$$W^{i_1\ldots i_k}_{~~i_1\ldots i_k} =
{W_{{}^{_{\!(q{=1})}}}}^{\!\!\!i_1\ldots i_k}_{\,i_1\ldots i_k}$$
since ${\Lambda_{{}^{_{\!(q{=1})}}}}^{\!\!\!ab}_{cd}$ never enters
${W_{{}^{_{\!(q{=1})}}}}^{\!\!\!i_1\ldots i_k}_{\,i_1\ldots i_k}$.
We also have $\forall\; i_1,\ldots i_N,\,j_1,\ldots j_N$
\begin{equation}
\epsilon_q^{i_1\ldots i_N}W^{1\ldots N}_{~~j_1\ldots j_N}
=W^{i_1\ldots i_N}_{~~j_1\ldots j_N} \;,\label{eWW}
\end{equation}
indeed, applying both sides to $dx^{j_1}\otimes\ldots dx^{j_N}$,
we obtain  the identity $\epsilon_q^{i_1\ldots
i_N}dx^{1}\wedge\ldots dx^{N} =dx^{i_1}\wedge\ldots dx^{i_N}$. If
all $j$'s are different we have $W^{j_1\ldots
j_N}_{~~j_1\ldots j_N}=1$ and, using (\ref{eWW}), we also have
\begin{eqnarray}
\label{We}& &W^{i_1\ldots i_N}_{~~1\ldots N}
=\epsilon_q^{i_1\ldots i_N}~,\\ \label{Wem1}& &W^{1\ldots
N}_{~~j_1\ldots j_N}
%=(\epsilon_q^{j_1\ldots j_N})^{-1}
={{\epsilon_{{}_{^{\!\!(q^{\!-\!1\!})}}}}}{}^{\!j_1\ldots j_N}
\equiv{{\epsilon_{{}_{^{\!\!(q^{\!-\!1\!})}}}}}{}_{j_1\ldots j_N}
\end{eqnarray}
where in the last line we have used  that
$\epsilon_q^{j_1\ldots j_N}$ is just a monomial in the $q_{ab}$'s
and that the inverse is therefore the same monomial with
$q_{ab}\rightarrow q_{ab}^{-1}$. The definition of
$\epsilon_{{}_{^{\!\!(q^{\!-\!1\!})}}}$ with lower indices
is just to preserve the index structure of
$W^{1\ldots N}_{j_1\ldots j_N}$. Relation (\ref{Wem1})
[as well as (\ref{We})] holds also for arbitrary $j$'s;
indeed both the l.h.s. and the r.h.s. are zero unless all
$j$'s are different. From (\ref{We}) and (\ref{Wem1}) we see that in
$W_{1\ldots N}$ the upper indices are $q$-antisymmetric, while the
lower indices are $q^{-1}$-antisymmetric; this is actually true
[cf. (\ref{eem1=W})] for any $W_{1\ldots k}$. In other words, in
the expression $\al={1\over k!}\al_{i_1\ldots
i_k}dx^{i_1}\wedge\ldots dx^i_k$ we can consider $\al_{i_1\ldots
i_k}$ with $q^{-1}$-antisymmetrized indices.

\medskip
\begin{proposition}
\begin{equation}
\label{eem1=W}
\epsilon_q^{i_1\ldots i_k l_{k+1}\ldots l_N}
{{\epsilon_{{}_{^{\!\!(q^{\!-\!1\!})}}}}}{}_{j_1\ldots j_k
l_{k+1}\ldots l_N}
=\sma{\mbox{$(N-k)!$}}\,W^{i_1\ldots i_k}_{~~j_1\ldots j_k}
\end{equation}
\end{proposition}
{\sl Proof.} We use an induction procedure. Indeed relation
(\ref{eem1=W}) holds for $k=N$; we consider it holds for $k$ and
we show it holds for $k-1$. We have to prove that
\begin{equation}
%\sum_{i_k}
W^{i_1\ldots i_{k-1}i_k}_{~~j_1\ldots j_{k-1}i_{k}}=
\sma{\mbox{$(N-k+1)$}}\,W^{i_1\ldots i_{k-1}}_{~~j_1\ldots j_{k-1}}
\label{nk1}
\end{equation}
or, by applying both sides of (\ref{nk1}) to
$dx^{j_1}\otimes\ldots dx^{j_{k-1}}$ and using (\ref{wedge2}),
\begin{equation}
\label{nk1_equiv} {\cal I}^{i_1\ldots i_{k-1} i_k}_{~b_1\ldots
b_{k-1}i_k} dx^{b_1}\wedge\ldots
dx^{b_{k-1}}={\sma{\mbox{$(N-k+1)$}}^{\,}} dx^{i_1}\wedge\ldots
dx^{i_{k-1}}~.
\end{equation}
Now from $\Lambda^{pi_k}_{~q i_k}=\delta^p_q$ and
(\ref{wedge4}) we have
$
\Tr_k {\cal I}_{1\ldots k}= N^{{}^{{}^{\,}}}I  - {\cal I}_{1\ldots k-1},
$
where $\Tr_k$ means trace on the $k$--factor of $\Omega^{\otimes
k}(\Re_q^N)$. Relation (\ref{nk1_equiv}) is then  proven by observing that
$\Lambda^{i_{v-1}i_{v}}_{~b_{v-1}b_v} dx^{b_{v-1}}\!\wedge
dx^{b_v}=-dx^{i_{v-1}}\!\wedge dx^{i_{v}}$. \fidi

\medskip
Notice that since the epsilon tensor up to a sign is invariant
under cyclic permutations [recall (\ref{prodotto_dei_q})] we also
have
\begin{equation}
\label{eem1bis=W}
\epsilon_q^{l_{k+1}\ldots l_N i_{1}\ldots i_k}
{{\epsilon_{{}_{^{\!\!(q^{\!-\!1\!})}}}}}{}_{l_{k+1}\ldots l_N
j_{1}\ldots j_k}
=\sma{\mbox{$(N-k)!$}}\,W^{i_{1}\ldots i_{k}}_{~~j_{1}\ldots j_{k}}
\end{equation}
and therefore
\begin{equation}
W^{i_1i_2\ldots i_k}_{~~i_1j_2\ldots j_{k}}=
\sma{\mbox{$(N-k+1)$}}\,W^{i_2\ldots i_{k}}_{~~j_2\ldots j_{k}}
\label{nk1mirror}~.
\end{equation}
\sk
The metric on  $\Re_q^{N}$ induces the following
pairing\footnote{
The shell structure of this pairing is uniquely determined by requiring
compatibility with the wedge product, see (\ref{W_Hermitian}).
The sign $(-1)^{[{k\over 2}]}=(-1)^{k(k-1)\over
2}$ is introduced in order to obtain in the commutative limit
the standard metric on $k$-forms.} $\langle~\,,~\rangle:\Omega_1^{\otimes
k}(\Re_q^N)\otimes\Omega_1^{\otimes k}(\Re_q^N)\rightarrow \Re_q^{N}$
(the tensor product $\otimes$ is over $\Re_q^{N}$):
\begin{equation}
\label{pairing}
\langle dx^{i_1}\otimes\ldots dx^{i_k},dx^{j_k}\otimes\ldots
dx^{j_1}\rangle = {(-1)^{[{k\over 2}]}\over k!}^{\,}g^{i_1j_1}\ldots g^{i_kj_k}\;,
\end{equation}
extended to all $\Omega_1^{\otimes k}(\Re_q^N)$ by $\langle f \th,\th'
\,h\rangle=f\langle \th,\th'\rangle h$ where $f\,,h\in\Re_q^{N}$
and
$\th,\th'\in \Omega_1^{\otimes k}(\Re_q^N)$. It is also easy to see that
$\langle~,~\rangle$ is a bimodule pairing, {\it i.e.}
\begin{equation}
\langle \th^{\,}f,\th'\rangle=\langle \th,f\th'\rangle~.
\label{goodpairing}
\end{equation}
In order to study the pairing between forms we need the following
properties among metric and epsilon tensors
\begin{lemma}
\label{metricandepsilon}
\begin{eqnarray}
\label{l1}
& &\det{}_qg\equiv\epsilon_q^{i_1\ldots i_N}\,g_{1 i_1}\ldots g_{N i_N}=
\det g=(-1)^{[{N\over 2}]}\\ \label{l2}
& &\epsilon_q^{j_1\ldots j_N}\,g_{j_1 i_1}\ldots g_{j_N i_N}=
\epsilon_q^{i'_1\ldots i'_N}= \epsilon_q^{i_1\ldots i_N} \det g\\
\label{l3}
& &{{\epsilon_{{}_{^{\!\!(q^{\!-\!1\!})}}}}}{}_{j_1\ldots j_N}
=\epsilon_q^{j_N\ldots j_1} \det g
\end{eqnarray}
\end{lemma}
{\it Proof}. Relation (\ref{l1}) follows from
$g_{ij}=\delta_{ij'}$ and $\epsilon_q^{N,N-1\ldots 1}=
(-1)^{[{N\over 2}]}\, \epsilon_q^{1\,2\ldots N}$, a consequence of
(\ref{epsilonqr}) and (\ref{relations2}).

In order to prove (\ref{l2}) we observe that if for a given $N$-tuple
$(i_1,\ldots k,l,\ldots i_N)$ we have $\epsilon_q^{i_1\ldots k
l\ldots i_N}= \epsilon_{_{\mbox{\scriptsize{$q$}}}}^{i'_1\ldots
k'l'\ldots i_N}\det g$ then also $\epsilon_q^{i_1\ldots l k\ldots
i_N}= \epsilon_{_{\mbox{\scriptsize{$q$}}}}^{i'_1\ldots l'k'\ldots
i_N}\det g$. Since $\epsilon_q^{N,N-1\ldots 1}=
\epsilon_q^{1\,2\ldots N}\det g$, relation (\ref{l2}) can be
proven by iterating this procedure.

Relation (\ref{l3}) follows from $\epsilon_q^{1\,2\ldots N}
={{\epsilon_{{}_{^{\!\!(q^{\!-\!1\!})}}}}}{}_{1\,2\ldots N}=
{{\epsilon_{{}_{^{\!\!(q^{\!-\!1\!})}}}}}{}_{N,N-1\ldots 1} \det
g$ and an iteration argument similar to the previous one. \fidi

\medskip
In the following proposition
we describe the  coupling between forms.

\begin{proposition} The pairing $\langle~,~\rangle$ satisfies the
following property
\begin{equation}
\label{W_Hermitian} \langle dx^{a_1}\otimes\ldots
dx^{a_k},dx^{i_k}\wedge\ldots dx^{i_1}\rangle= \langle
dx^{a_1}\wedge\ldots dx^{a_k},dx^{i_k}\otimes\ldots
dx^{i_1}\rangle ~,
\end{equation}
and when restricted to forms reads
\begin{equation}
\langle dx^{a_1}\wedge\ldots dx^{a_k},dx^{i_k}\wedge\ldots
dx^{i_1}\rangle \label{pairing_on_R} =(-1)^{[{k\over 2}]}\,g^{a_k
b_k}\ldots g^{a_1 b_1} \,W^{i_k\ldots i_1}_{~~b_k\ldots b_1}~.
\end{equation}
\end{proposition}
{\sl Proof. } Relation (\ref{W_Hermitian}) is equivalent
to
\begin{equation}
\label{W_hermitian_equiv} g^{a_k b_k}\ldots g^{a_1 b_1}
\,W^{i_k\ldots i_1}_{~~b_k\ldots b_1} = W^{a_1\ldots
a_k}_{~~b_1\ldots b_k} \, g^{b_1 i_1}\ldots g^{b_k i_k} \;,
\end{equation}
which holds for $k=N$ because of (\ref{eem1=W}), and the previous
lemma. Now, by induction, if (\ref{W_hermitian_equiv})
holds for $k$ then it holds for
$k-1$: indeed just multiply by $g_{a_ki_k}$ (summing over both
$a_k$ and $i_k$) and recall (\ref{nk1}) and (\ref{nk1mirror}).
\fidi

\medskip
The $*$--Hodge operator is defined as the unique map
$*:\Omega_k(\Re_q^{N})\rightarrow\Omega_{N-k}(\Re_q^{N})$
such that
\begin{equation}
\label{defhodge}
\alpha\wedge*\beta = \langle\alpha,\beta\rangle V_{N}
~~~~~~\alpha,\beta\in\Omega_k(\Re_q^{N})\;.
\end{equation}
We collect in the following proposition the main  properties of the
$*$--Hodge operator.

\begin{proposition}
\label{Hodge} An explicit expression for $*$ is given by
\begin{equation}
\label{h1}*(dx^{i_1}\wedge\ldots dx^{i_k})= C_{N,k} \;
\epsilon_q^{i_1\ldots i_k l_{k+1}\ldots l_N}
g_{l_{k+1}t_{k+1}}\ldots g_{l_{N}t_{N}}\, dx^{t_N}\wedge\ldots
dx^{t_{k+1}},\end{equation} where $C_{N,k}=(-i)^{[{N\over
2}]}\,(-1)^{[{N-k\over 2}]}/ (N-k)!$.

{}Forall $\al,\be\in\Omega_k(\Re_q^N)$, $f,h\in\Re_q^N$ and
$\ga\in\Omega_{N-k}(\Re_q^N)$ we have
\begin{eqnarray}
& &\nonumber\\[-1em]
\label{h2}& &*(f\al\, h)=f(*\al)h~~~~~~~~~~~\sma{$\Re_q^{N}$ left and right linearity,}\\
%bimodule ${}^{\!}$map property,}\\
\label{h3}& &*1=V_N~~,~~~~*V_N=1~~,\\
\label{h4}& & **\al=(-1)^{k(N-k)}\,\al~~,\\
\label{h5}& &\al\wedge *\be =(-1)^{k(N-k)}\,*\al\wedge\be~~,\\
\label{h6}& &\langle\alpha,\beta\rangle = \langle *\alpha,*\beta\rangle~~,\\
\label{h7}& & \langle *{}^{\!}\,\alpha,\ga\rangle = \langle\alpha\wedge\ga,V_N\rangle~.
\end{eqnarray}
\end{proposition}

\noi {\sl Proof. }
In order to prove (\ref{h1}) we have to show that
$dx^{a_1}\wedge\ldots dx^{a_k}\wedge *(dx^{i_1}\wedge\ldots dx^{i_k})
= \langle dx^{a_1}\wedge\ldots dx^{a_k},
dx^{i_1}\wedge\ldots dx^{i_k}\rangle V_{N}$
where in the l.h.s. we use (\ref{h1}). Equivalently
we have to show
\begin{eqnarray}
 {{(-1)^{[{N\over 2}]}}\,(-1)^{[{N-k\over 2}]}\over (N-k)!}\,
\epsilon_q^{i_1\ldots i_k l_{k+1}\ldots l_N}
g_{l_{k+1}t_{k+1}}\ldots g_{l_{N}t_{N}}\, \epsilon_q^{{a_1}\ldots
{a_k}{t_N}\ldots {t_{k+1}}}=
&\nonumber\\
&\!\!\!\!\!\!\!\!\!\!\!\!\!\!\!\!\!\!\!\!\!\!\!\!\!\!\!\!\!{}(-1)^{[{k\over 2}]}
g^{a_kb_k}\ldots g^{a_1b_1} \,W^{i_1\ldots i_k}_{~~b_k\ldots b_1}
\;. \nonumber
\end{eqnarray}
This last equality can be proven multiplying by $g_{c_1a_1}\ldots
g_{c_ka_k}$, using (\ref{l2}) and (\ref{l3}), recalling that the
$q$-epsilon tensor is  invariant (up to a sign) under cyclic
permutations and finally using (\ref{eem1=W}). (\ref{h2}) follows
from (\ref{goodpairing}). The second relation in (\ref{h3})
follows from $\langle V_N,V_N\rangle=1$. (\ref{h4}) can be proven
as in the commutative case, for example using twice (\ref{h1}) and
then (\ref{l2}), (\ref{l3}) and $W_{1\ldots k}W_{1\ldots k}=k!
W_{1\ldots k}$. Also (\ref{h5}) can be proved using (\ref{h1}) as
in the commutative case. (\ref{h6}) and (\ref{h7}) are then easily
shown using (\ref{defhodge}) and (\ref{h5}). \fidi

\medskip
\begin{remark}{\rm It is easy to verify that the $*$--Hodge operator
defined in (\ref{defhodge}) coincides with the one
defined in \cite{ConnesDuboisViolette}.
% induced by the
%$*_o$ defined on the commutative plane via the splitting
%homomorphism in \cite{ConnesDuboisViolette}.
Indeed let's denote
the relations on the $n$--torus $T^n_\theta$ by
$U^iU^j=q_{ij}U^jU^i$, ${U^i}^*=U^{i'}={U^i}^{-1}$ and $U^i=1$ if
$i=i'$, where $i=1,\ldots N,\; N=2n$ or $N=2n+1$ and $q_{\al\be}=
e^{i\theta_{\al\be}}$
$\al,\be=1,\ldots n$. The exterior algebra $\Omega(\Re_q^N)$
is then homomorphic to  $\Omega(\Re^N)\otimes  C_{alg}(T^n_\theta)$
via the identification
$dx^{i_1}\wedge\ldots dx^{i_k}=dx^{i_1}_0\wedge\ldots dx^{i_k}_0\otimes
U^{i_1}\ldots U^{i_k}$ where the $x_0$'s are the coordinates on the
commutative plane. In particular [cf. (\ref{epsilonqr})] we have
$\epsilon_q^{i_1\ldots i_N}=\epsilon^{i_1\ldots i_N}U^{i_1}\ldots
U^{i_N}$.
By applying $*_0\otimes id$ to $dx^{i_1}_0\ldots dx^{i_k}_0\otimes
U^{i_1}\ldots U^{i_k}$ (where $*_0$ is the commutative $*$--Hodge operator)
we obtain expression (\ref{h1}) thus showing that $*=*_0\otimes\id$. $~$\fidi}
\end{remark}
\medskip
\noi{\bf Conjugation}\\ We now study the star structure that
$\Omega(\Re_q^N)$ inherits from $\Re_q^N$. We recall that
$(x^a)^\*=x^{a'}$ and that on $1$-forms $(fdh)^\*=d(h^\*)\;f^*$.
On $\Omega_1^{\otimes k}(\Re_q^N)$ we define \cite{Woronowicz}
$$(dx^{i_1}\otimes\ldots dx^{i_k})^\*\equiv (-1)^{k(k-1)\over 2}
dx^{i'_k}\otimes\ldots dx^{i'_1}~.$$
Notice that the complex conjugate of $\epsilon_q^{i_i\ldots i_N}$
is ${{\epsilon_{{}_{^{\!\!(q^{\!-\!1\!})}}}}}{}_{i_1\ldots i_N}$,
so that using (\ref{l3}) and (\ref{l2})
${\overline{\epsilon_q^{i_1\ldots i_N}}}=\epsilon_q^{{i'}_N\ldots {i'}_1}$,
then from (\ref{eem1=W}) and (\ref{eem1bis=W}) we have
${\overline{W^{i_{1}\ldots i_{k}}_{~~j_{1}\ldots j_{k}}}}=
W^{{i'}_{k}\ldots {i'}_{1}}_{~~{j'}_{k}\ldots {j'}_{1}}$ {\it i.e.}
\begin{equation}
(dx^{i_1}\wedge\ldots dx^{i_k})^\*= (-1)^{k(k-1)\over 2}
dx^{i'_k}\wedge\ldots dx^{i'_1}~,\label{ultima}
\end{equation}
 thus the space of exterior forms $\Omega(\Re_q^N)$
 naturally inherits the conjugation on $\Omega_1^{\otimes k}(\Re_q^N)$.
 In particular we have that the volume element $V_N$ is real:
 ${V_N}^\*=V_N$.

\begin{proposition}Reality of the $*$--Hodge operator. For any form $\al$
\begin{equation}\label{h8}*\;\al^\*=(*\al)^\*
\end{equation}
\end{proposition}
{\sl Proof.}
(\ref{h8}) can be proven by explicit computation, using (\ref{h1}),
recalling (\ref{l2}),(\ref{l3}) and again
the invariance (up to a sign) of the
$q$-epsilon tensor under cyclic permutations.
\fidi

\medskip

\subsection{Hodge Theory on $\Sq^N$}
%, in a chart $U\subset \S^N$, $\Omega(\S^N)\subset
%\Omega\left.(\Re^{N+1})\right|_{{\S^{{}^{_N}}}}$.
Recalling that the versor normal to the commutative
unit sphere $\S^N$ is $dc/2$,
it is easy to see that forall $\th,\th'\in \Omega(\S^N)$, the metric
on $\Omega(\S^N)$ and on $\Omega(\Re^{N+1})$ are related by
$\left.\langle \th,\th'\rangle\right|_{{\S^{{}^{_N}}}}
=\left.\langle \th\wedge {dc\over 2},\th'
\wedge {dc\over 2} \rangle\right|_{{\S^{{}^{_N}}}}$.
It is therefore natural to define in the noncommutative case,
forall $[\al],[\be]\in\Omega_k(\Sq^N)$,
\begin{equation}
\label{pairing_on_S} \langle [\al],[\be]\rangle = {1\over 4}[\langle
\al\wedge dc,\be\wedge dc \rangle]~~.
\end{equation}
Independence from the representatives $\al$ and  $\be$ is easily proven.
The $*$--Hodge map is then the unique map
$*:\Omega_k(\Sq^{N})\rightarrow\Omega_{N-k}(\Sq^{N})$
such that
\begin{equation}
\label{shodgedef}
{}~~~~[\alpha]\wedge*[\beta] = \langle[\alpha],[\beta]\rangle \,\V_N
~~~~~~~[\alpha],[\beta]\in\Omega_k(\Sq^{N})\;.
\end{equation}
It is also easy to check that  $J$
is a $\*$-ideal: $J^\*\subset J$. Then
\begin{equation}
[\al]^\*=[\al^\*]
\end{equation}
is a well defined $\*$-structure on $\Omega(\Sq^{N})$. Since $J$ is a
differential ideal and a left coideal we have
$d\,[\al]^\*=(d[\al])^\*$ and we have a $SO_q(N,\Re)$-coaction
on $\Omega(\Sq^{N})$.
Reality of the volume form $\V_N$ follows from reality of
$V_{N+1}$ and of $c^{\,}$; we have
$$\begin{array}{c}
x^k\omega_k\wedge dc=2cV_{N+1}=(2cV_{N+1})^\*
=(-1)^Ndc^\*\wedge\omega_k^{\,\*}{x^k}^\*
=\omega_k^{\,\*}{x^k}^\*\wedge dc
\end{array}$$
and therefore
$
\,\V_N^{\,\*}=[x^k\omega_k]^\*=[\omega_k^{\,\*}{x^k}^\*]=[x^k\omega_k]=
\V_N~.
$

We collect in the following proposition the main properties of the
$*$--Hodge operator on $\Omega(\Sq^{N})$.

\begin{proposition}
\label{sHodge} An explicit expression for $*$ is given by
\begin{equation}
\label{sh1}*[dx^{i_1}\wedge\ldots dx^{i_k}]=
C'_{N,k} \;
\epsilon_q^{i_1\ldots i_k a l_{k+1}\ldots l_N}
g_{ab}^{\,}g_{l_{k+1}t_{k+1}}\ldots g_{l_{N}t_{N}}\, [dx^{t_N}\wedge\ldots
dx^{t_{k+1}}\,x^b]\,,
\end{equation}
where $C'_{N,k}=(-i)^{[{N+1\over
2}]}\,(-1)^{[{N-k\over 2}]}(-1)^{N-k}/ (N-k)!$.

{}Forall $\be\in\Omega_k(\Re_q^{N+1})$,
$\th,\eta\in\Omega_k(\Sq^N)$, $f,h\in\S_q^N$ and
$\nu\in\Omega_{N-k}(\S_q^N)$ we have
\begin{eqnarray}
& &\nonumber\\[-1em]
\label{sh1'}& &*[\be]=(-1)^{N-k}[*(\be\wedge {dc\over 2})]~~,\\
%~~~~~~~~~~~\sma{where $\al\in\Omega_k(\Re_q^N)$ and $[\al]=\th\,$,}\\
\label{sh2}& &*(f\th\, h)=f(*\th)h~~~~~~~~~~~~~~~~~~\sma{$\Sq^{N}$
left and right linearity}\,,\\
\label{sh3}& &*1=\V_N~~,~~~~*\V_N=1~~,\\
\label{sh4}& & **\th=(-1)^{k(N-k)}\,\th~~,\\
\label{sh5}& &\th\wedge *\eta =(-1)^{k(N-k)}\,*\th\wedge\eta~~,\\
\label{sh6}& &\langle\th,\eta\rangle = \langle *\th,*\eta\rangle~~,\\
\label{sh7}& & \langle *{}^{\!}\,\th,\nu\rangle
= \langle\th\wedge\nu,\V_N\rangle~,\\
\label{sh8}& &*\;\th^\*=(*\th)^\*~~~~~~~~~~~~~~~~~~~~~~~~\sma{reality of the $*$-Hodge}\,.
\end{eqnarray}
\end{proposition}
{\sl Proof. } Relation (\ref{sh1'}) is equivalent to
$(-1)^{N-k}[\al]\wedge[*(\be\wedge dc/2)]=\langle [\al],[\be]\rangle \V_N$
{\it i.e.} (use $[c]=1$)
$$
\begin{array}{c}
{1\over 2}^{\,}(-1)^{N-k}[c^{\,}\al\wedge *(\be\wedge dc)]=
{1\over 4}[\langle \al\wedge dc,\be\wedge dc\rangle x^k\om_k]~;
\end{array}
$$
this last relation holds because
$$
\begin{array}{l}
\left\{{1\over 2}^{\,}(-1)^{N-k}{}^{\,}{c}^{\,}\al\wedge *(\be\wedge dc)
-{1\over 4}\langle \al\wedge dc,\be\wedge dc\rangle x^k\om_k
\right\}\wedge dc~=\\[.5em]
{}~~~~~~~~~~~~~~~~~~~~~~~~~~~~~~~~~~~~~~{1\over 2}^{\,}{c}{\,}
\al\wedge dc\wedge *(\be\wedge dc)
-{1\over 2} \langle \al\wedge dc,\be\wedge dc\rangle^{\,} c^{\,} V_{N+1}
=0~.
\end{array}
$$
Relation (\ref{sh1}) follows from (\ref{sh1'}) and (\ref{h1}).
Relation (\ref{sh2}) follows from (\ref{shodgedef}) and
(\ref{goodpairing}). Also (\ref{sh3}) easily follows from
 (\ref{shodgedef}).
Relation (\ref{sh4}) is equivalent to
\begin{equation}
{1\over 4}(-1)^{N-k}\,(-1)^k[*(*(dx^{i_1}\wedge\ldots dx^{i_k}\wedge dc)\wedge dc)]
=(-1)^{k(N-k)}[dx^{i_1}\wedge\ldots dx^{i_k}]~.\label{eqc}
\end{equation}
Applying twice (\ref{h1}) and then (\ref{l2}) and (\ref{l3})
the l.h.s. equals
$$
\begin{array}{ll}
{(-1)^{k(N-k)}}^{\,}{1\over k!}^{\,}
W^{i_1\ldots i_k a}_{~~s_1\ldots s_k f}
[dx^{s_1}\wedge\ldots dx^{s_k}g_{ab}x^fx^b]=&\\[.5em]
&\!\!\!\!\!\!\!\!\!\!\!\!\!\!\!\!\!\!\!\!\!\!\!\!\!\!\!\!
\!\!\!\!\!\!\!\!\!\!\!\!\!\!\!\!\!\!\!\!\!\!\!\!\!\!\!\!
\!\!\!\!\!\!\!\!\!\!\!\!\!\!\!\!\!\!\!\!\!\!\!\!\!\!\!\!
(-1)^{k(N-k)}{}_{\,}[dx^{i_1}\wedge\ldots dx^{i_k}{}^{\,}c
-{\cal I}^{i_1\ldots i_{k-1}i_k}_{~~u_1\ldots u_{k-1}f}{}_{\,}
dx^{u_1}\wedge\ldots dx^{u_{k-1}}\wedge dc\;x^f]
\end{array}$$
where in the second line we used (\ref{wedge2}),(\ref{wedge4}) and
$\Lambda^{va}_{~~u_k f}=q_{va}\delta^v{}_f\delta^a{}_{u_k}$. Finally
notice that this last expression is  the same equivalence class
as the one in the r.h.s. of (\ref{eqc}).
Relation (\ref{sh5}) is equivalent to
$[\al\wedge *(\be\wedge dc)]=(-1)^{k(N-k)}
[*(\al\wedge dc)\wedge \be]$ (with $[\al]=\th,$  $[\be]=\eta$).
This equality holds if
$$\begin{array}{c}
\al\wedge *(\be\wedge dc)\wedge dc=(-1)^{k(N-k)}
*(\al\wedge dc)\wedge \be\wedge dc
\end{array}
$$
{\vskip -0.6cm}
{\it i.e.}
{\vskip -0.6em}
$$\begin{array}{c}
(\al\wedge dc)\wedge *(\be\wedge dc)=(-1)^{(k+1)(N-k)}
*{\!}(\al\wedge dc)\wedge (\be\wedge dc)
\end{array}
$$
and this last relation is true because of (\ref{h5}) (with
\sma{$N+1$} instead of \sma{$N$}).
(\ref{sh6}) and (\ref{sh7}) are then easily
shown using (\ref{shodgedef}) and (\ref{sh5}).
The reality property of the $*$--Hodge operator on $\Omega(\Sq^N)$
follows from that on $\Omega(\Re_q^{N+1})$ and
from $(\al\wedge dc)^\*=\al^\*\wedge dc$.
\fidi

\bigskip
\section{Geometry of the instanton bundle on $\Sq^{2n}$}
The purpose of this section is to apply the calculus developed in
the previous sections to the twisted sphere $\Sq^{2n}$, namely
to study the geometric properties of the instanton projector
introduced in \cite{ConnesLandi} for the $4$--dimensional case and in
\cite{ConnesDuboisViolette} for the general case. As shown in
\cite{ConnesDuboisViolette} its curvature is self dual -- for a different
and more explicit proof in the $4$--dimensional case see the first
version of this paper, math.QA/0108136v1. We here
compute the charge, {\it i.e.} the
Chern--Connes pairing between the projector $e$ and the character
$\tau$ defined in (\ref{carattere}).

We briefly recall some basic notions about the coupling between
cyclic homology and $K$--theory, (general references for this
section are \cite{ConnesBook} and \cite{Loday}). Let $A$ be an
associative $\Co$--algebra. A projector $e\in M_k(A)$, {\it i.e.}
$e^2=e$, defines a finitely generated projective module ${\cal E}
= e A^k$, which we consider as the space of sections of a {\it
quantum vector bundle}. In the Grothendieck group $K_0(A)$ the
module $\cal E$ defines a class, that we still denote with $e$.
Topological informations of this quantum bundle are extracted by
using cyclic homology and cohomology, which are defined as
follows. Let $d_i:A^{\otimes (n+1)}\rightarrow A^{\otimes n}$ be
defined as $d_i(a_0\otimes a_1\ldots \otimes a_n
)=a_0\otimes..\,a_ia_{i+1}.. \otimes a_n$, for $i=0,\ldots n-1$
and $d_n(a_0\otimes a_1\ldots \otimes a_n)=a_na_0\otimes a_1\ldots
\otimes a_{n-1}$; the Hochschild boundary is defined as
$b=\sum_{i=0}^n(-1)^id_i$ and the Hochschild complex is
($C_*(A),b$), with $C_n(A)=A^{\otimes n+1}$.  Let $t_n(a_0\otimes a_1\otimes
\ldots a_n)=(-1)^n a_1\otimes\ldots a_n\otimes a_0$ be the cyclic
operator and $C^\lambda_n(A)=A^{\otimes n+1}/(1-t)A^{\otimes
n+1}$. The Connes complex is then ($C^\lambda_*(A),b$); its
homology is denoted as $H^\lambda_*(A)$. Analogously we define
$H_\lambda^*(A)$ as the cohomology of the complex
($C_\lambda^*(A),b$), where
$C^n_\lambda=\{\tau:A^{\otimes(n+1)}\rightarrow\Co\,|\, \tau\sma{$^{\,}\circ^{\,}$}
t_n=(-1)^n\tau\}$ and $b(\tau)=\tau\sma{$^{\,}\circ^{\,}$} b$.

Cycles give an alternative way of introducing cyclic cocycles. A
cycle ($\Omega,d,\int$) over $A$ of dimension $n$ is given by a
differential graded algebra
($\Omega=\bigoplus_{k=0}^n\Omega_k,d$), a closed graded trace
$\int:\Omega_n\rightarrow \Co$ and an algebra morphism
$\rho:A\rightarrow\Omega_0$. More explicitly these data must
satisfy:
\begin{eqnarray*}
d(\omega\nu)=(d\omega) \nu + (-1)^{|\omega|}\omega d\nu ~,~~~~~~~
d^2=0\cr \int\omega\nu=(-1)^{|\omega||\nu|}\int\nu\omega ~,~~~~~~~
\int d\omega=0  \;.
\end{eqnarray*}
The character of the cycle ($\Omega,d,\int$) over $A$, defined as
$\tau:A^{\otimes (n+1)}\rightarrow \Co\;,$
$$\tau(a_0\otimes\ldots  a_n)=\int\rho(a_0)d\rho(a_1)\ldots
d\rho(a_n)\;,$$ is a cyclic $n$--cocycle; we still denote its class in
$H^n_\lambda$ with $\tau$.
It can be shown that all
cyclic cocycles are characters of some cycle.

For each projector $e\in M_k(A)$, {\it i.e.} $e^2=e$, the Chern
character is a class in cyclic homology and is defined as
$ch_n^\lambda(e)=1/n!^{\,}\Tr(e^{\otimes 2n+1})\in
H^\lambda_{2n}(A)$, where $\Tr:M_k(A)^{\,\otimes n}\to A^{\otimes
n}$ is the generalized trace, {\it i.e.} $\Tr(N^{(1)}\otimes \ldots
\otimes N^{(n)})=
\sum_{\alpha}N^{(1)}_{\alpha_1\alpha_2}\otimes N^{(2)}_{\alpha_2\alpha_3}
\ldots
\otimes N^{(n)}_{\alpha_n\alpha_1}$. For each $2n$--cocycle $\tau$ the
Chern--Connes pairing given by
$\langle e,\tau \rangle = \tau(ch_n^\lambda(e))$ depends only
on the class of $e$ in $K_0(A)$ and of $\tau$ in $H^{2n}_\lambda$ and
can be computed in the following way. Let ($\Omega,d,\int$) be the cycle over $A$ of dimension $2n$
associated to $\tau$. We canonically define a connection on ${\cal
E}=e A^k$ in the following way. Let $f_\alpha=e v_\alpha$, where
$(v_\alpha)_\beta=\delta_{\alpha\beta}$; the Levi--Civita
connection $\nabla^e:
{\cal E}\rightarrow{\cal E}\otimes_A\Omega_1$ is defined as
$$
\nabla^e(f_\alpha) = \sum_\beta f_\beta \otimes de_{\beta\alpha} \;.
$$
The curvature ${\cal F}^e=\nabla^2:{\cal E}\rightarrow {\cal
E}\otimes \Omega_2$ is  $A$--linear and defines a $k\times
k$--matrix of two forms given by ${\cal F}^e=edede$. The
Chern--Connes pairing between the character of the cycle and
${\cal E}$ is computed by the following formula:
\begin{equation}
\label{carattere_col_ciclo} \langle e,\tau \rangle =
{1\over n!}^{\,}\tau(\Tr(e^{\otimes 2n+1}))=\frac{1}{n!} \int \Tr(({\cal F}^e)^n) \;.
\end{equation}

\medskip
Let us come back to $\Sq^{2n}$ and let us introduce the
Clifford algebra defined in \cite{ConnesDuboisViolette}. We have
that ${\rm Cliff}(\Re_q^{2n+1})$ is generated
by $2n+1$ generators $\gamma^i$ and the following relations:
\begin{equation}
\label{clifford}
\gamma^i \gamma^j + q_{ji} \gamma^j \gamma^i = 2 g^{ij}\;;
\end{equation}
let us remark that the chiral $\gamma$ of \cite{ConnesDuboisViolette}
is included among the fundamental $\gamma$'s and corresponds to
$\gamma^{n+1}$. The unique irreducible representation is given on
$\otimes^n\Co^2$ by
\begin{equation}\label{clifford_rep}\gamma^i = \sqrt{2}
\left(\begin{array}{cc}-q_{i1}&0\cr0&1\end{array}\right)
\otimes\ldots\otimes
\left(\begin{array}{cc}-q_{i,i-1}&0\cr0&1\end{array}\right)\otimes
\left(\begin{array}{cc}0&0\cr 1&0\end{array}\right)\otimes 1\ldots 1\;,
~~~~~ i<n
\end{equation}
$$\gamma^{i'} = \gamma^i{}^\dagger ~~,~~~~~~~ \gamma^{n+1} = \otimes^n
\left(\begin{array}{cc}1&0\cr0&-1\end{array}\right) \;.
$$
As a consequence of (\ref{clifford}) $e = \frac{1}{2}
\left(1+\gamma^i [x^{j}]g_{ij}\right)$ is a projector.

Let ($\Omega(\S^{2n}_q),d,\int$) be the cycle over $\S^{2n}_q$ defined
in Proposition \ref{cycle}, and let $\tau$ be its character
defined in (\ref{carattere}). Let $\nabla^e$ be the canonical
connection with values in $\Omega(\Sq^{2n})$ and let ${\cal
F}^e= edede$ be its curvature. As a consequence of $e=e^\*$
and of (\ref{ultima}) ${\cal F}^e$ is antihermitian, {\it i.e.}
${{\cal F}^{e}_{\alpha\beta}}^\*=-{\cal F}^{e}_{\beta\alpha}$,
$\,\alpha,\beta=1,\ldots 2^n$.
Moreover in \cite{ConnesLandi} and in
\cite{ConnesDuboisViolette} it is shown that
$ch_m^\lambda(e)=0$ for $m<n$. In the
following proposition we compute
the maximal Chern character and verify that the normalization of
$\tau$ discussed in (\ref{carattere}) still
guarantees the integrality of the pairing.

\medskip
\begin{proposition}
The charge of the instanton projector on $\sphere_q^{2n}$,
{\it i.e.} the Chern--Connes pairing between $e$ and $\tau$, reads
$$\langle e,\tau\rangle = \frac{1}{n!}\ \tau(\Tr[e^{\otimes 2n+1}]) = 1\;.$$
\end{proposition}
{\sl Proof.} By the use of
the faithful representation (\ref{clifford_rep}) it can be shown that
$$
\Tr (\gamma_{i_0}\ldots\gamma_{i_{2n}}) =
2^n {{\epsilon_{{}_{^{\!\!(q^{\!-\!1\!})}}}}}{}_{i_0\ldots i_{2n}}\;.
$$
Since thanks to Stokes' Theorem we can ignore $\Tr([de]^{2n})$
(it can be shown that it is zero anyway), we
have that
\begin{eqnarray*}
\int\Tr [e(de)^{2n}] &=& \int\frac{1}{2^{2n+1}}
\Tr(\gamma_{i_0}\ldots \gamma_{i_{2n}}) [x^{i_0'}dx^{i_1'}\ldots
dx^{i_{2n}'}] \cr
&=&\int\frac{1}{2^{n+1}}\
{{\epsilon_{{}_{^{\!\!(q^{\!-\!1\!})}}}}}{}_{i_0\ldots i_{2n}}
[x^{i_0'}dx^{i_1'}\ldots dx^{i_{2n}'}] \cr
&=&\int\frac{(-1)^n}{2^{n+1}}\
{{\epsilon_{{}_{^{\!\!(q^{\!-\!1\!})}}}}}{}_{i_0\ldots i_{2n}}
[x^{i_0}dx^{i_1}\ldots dx^{i_{2n}}] \cr
&=&\int \frac{(2n)!}{2^{n+1}}\ i^n [x^i\omega_i] =
\int \frac{(2n)!}{2^{n+1}}\ i^n \V_{2n}  \;,
\end{eqnarray*}
where we used (\ref{l2}) in the third line and the definition
(\ref{omega}) in the last line. The result then follows by applying the definition
(\ref{carattere}) of $\tau$. \fidi

\sk
\sk
\noi{\bf Acknowledgements}

\noi This research has been in part supported by a Marie Curie Fellowship
of the European Community programme IHP under contract number MCFI-2000-01982.


\begin{thebibliography}{99}

\bibitem{AschieriCastellani}
P.~Aschieri and L.~Castellani, ``An Introduction to noncommutative
differential geometry on quantum groups,'' Int.\ J.\ Mod.\ Phys.\
A {\bf 8} (1993) 1667 [arXiv:hep-th/9207084].
%%CITATION = HEP-TH 9207084;%%

%\cite{Aschieri:1996th}
\bibitem{Aschieri}
P.~Aschieri and L.~Castellani,
``Bicovariant Calculus on Twisted ISO(N), Quantum Poincare' Group and Quantum Minkowski Space,''
Int.\ J.\ Mod.\ Phys.\ A {\bf 11} (1996) 4513
[q-alg/9601006].
%%CITATION = Q-ALG 9601006;%%

%\cite{Aschieri:1999et}
\bibitem{Scarfone}
P.~Aschieri, L.~Castellani and A.~M.~Scarfone,
``Quantum orthogonal planes: $ISO_{q,r} (n+1, n-1)$ and $SO_{q,r} (n+1, n-1)$ bicovariant calculi,''
Eur.\ Phys.\ J.\ C {\bf 7} (1999) 159
[q-alg/9709032].
%%CITATION = Q-ALG 9709032;%%


\bibitem{NOI}
%\cite{Bonechi:2000ia}
%\bibitem{Bonechi:2000ia}
F.~Bonechi, N.~Ciccoli and M.~Tarlini,
``Noncommutative instantons on the 4-sphere from quantum groups,''
math.qa/0012236.
%%CITATION = MATH.QA 0012236;%%



\bibitem{Carow-Watamura:1991zp}
U.~Carow-Watamura, M.~Schlieker and S.~Watamura,
``SO-q(N) covariant differential calculus on quantum space and
quantum deformation of Schr\"odinger equation,''
Z.\ Phys.\ C {\bf 49} (1991) 439.
%%CITATION = ZEPYA,C49,439;%%


\bibitem{ConnesBook} A.~Connes, ``Noncommutative Geometry''.
Academic Press, (1994).


\bibitem{ConnesDuboisViolette} A.~Connes, M.~Dubois--Violette,
``Noncommutative finite-dimensional manifolds. I. Spherical
manifolds and related examples,'' math.QA/0107070.


%\cite{Connes:2001tj}
\bibitem{ConnesLandi}
A.~Connes and G.~Landi,
``Noncommutative manifolds: The instanton algebra and isospectral  deformations,''
Commun.\ Math.\ Phys.\  {\bf 221} (2001) 141
[math.qa/0011194].
%%CITATION = MATH.QA 0011194;%%


%\cite{Dobrev}
\bibitem{Dobrev}
V.~K.~Dobrev,
``Canonical q deformations of noncompact Lie (super)algebras,''
J.\ Phys.\ AA {\bf 26} (1993) 1317.
%%CITATION = JPAGB,A26,1317;%%


\bibitem{FRT}
L.~D.~Faddeev, N.~Y.~Reshetikhin and L.~A.~Takhtajan,
``Quantization Of Lie Groups And Lie Algebras,''
Lengingrad Math.\ J.\  {\bf 1} (1990) 193
[Alg.\ Anal.\  {\bf 1} (1990) 178].
%%CITATION = 00064,1,193;%%



\bibitem{Fiore:1994vk}
G.~Fiore,``Quantum groups $SO_q(N)$, $Sp_q(n)$ have $q$-determinants, too''
J.\ Phys.\ AA {\bf 27} (1994) 3795.
%%CITATION = JPAGB,A27,3795;%%


\bibitem{VarillyBook}
%\cite{Gracia-Bondia:2001tr}
%\bibitem{Gracia-Bondia:2001tr}
J.~M.~Gracia-Bondia, J.~C.~Varilly and H.~Figueroa,
``Elements of noncommutative geometry''.
Birkhaeuser, Boston, (2001).


\bibitem{JurcoWat} B.Jurco ``Differential Calculus on Quantized
Simple Lie Groups'' Lett.Math.Phys. {\bf 22} 177 (1991).

U.Carow--Watamura, M.Schlieker, S.Watamura, W.Weich, ``Bicovariant
and Differential Calculus on Quantum Groups $SU_q(n)$ and
$SO_q(n)$'' Commun.Math.Phys. {\bf 142} 605 (1991)

B.~Jurco,
``Differential calculus on quantum groups: Constructive procedure,''
in Proceedings of the International School of Physics ``Enrico Fermi'',
Varenna 1994, Course CXXVII: Quantum Groups and their Application in
Physics, L. Castellani, J. Wess eds. (IOS/Ohmsha Press 1996)
[hep-th/9408179].
%%CITATION = HEP-TH 9408179;%%

\bibitem{Loday} J.~L.~Loday,``Cyclic Homology''. Springer--Verlag
Berlin, (1992).


\bibitem{MNW} T.~Masuda, Y.~Nagakami, J.~Watanabe,
``Noncommutative Differential Geometry on the Quantum Two Sphere of
Podle\`s. I; An Algebraic Viewpoint,'' K-Theory {\bf 5} (1991) 151.



\bibitem{Ogievetsky:1992qp}
O.~Ogievetsky and B.~Zumino,
``Reality in the differential calculus on q Euclidean spaces,''
Lett.\ Math.\ Phys.\  {\bf 25} (1992) 121
[hep-th/9205003].
%%CITATION = HEP-TH 9205003;%%

%\cite{Podles:1987wd}
\bibitem{Podles:1987wd}
P.~Podle\`s,
``Quantum Spheres,''
Lett.\ Math.\ Phys.\  {\bf 14} (1987) 193.
%%CITATION = LMPHD,14,193;%%

%\cite{Wess:1991vh}
\bibitem{Wess:1991vh}
W.~Pusz and S.~L.~~Woronowicz, ``Twisted Second quantization,''
Rep. Math. Phys. {\bf 27} (1990) 231.

J.~Wess and B.~Zumino,
``Covariant Differential Calculus on the Quantum Hyperplane,''
Nucl.\ Phys.\ Proc.\ Suppl.\  {\bf 18B} (1991) 302.
%%CITATION = NUPHZ,18B,302;%%


\bibitem{Schirrmacher}
N. Reshetikhin,
``Multiparametric
quantum groups and twisted quasitriangular Hopf algebras,'' Lett. Math. Phys.
{\bf 20} (1990) 331.

%\cite{Schirrmacher:1991xx}
%\bibitem{Schirrmacher:1991xx}
A.~Schirrmacher,
``Multiparameter R matrices and its quantum groups,''
J.\ Phys.\ AA {\bf 24} (1991) L1249.
%%CITATION = JPAGB,A24,L1249;%%



\bibitem{Sitarz} A.~Sitarz,  ``Twists and spectral triples for
isospectral deformations,'' math.QA/0102074

%\cite{Steinacker}
\bibitem{Steinacker}
H.~Steinacker,
``Integration on quantum Euclidean space and sphere in N-dimensions,''
J. Math. Phys. {\bf 37} (1996) 7438 [q-alg/9506020].


\bibitem{Twietmeyer:1992mj}
E.~Twietmeyer,
``Real forms of $U_q(g)$,''
Lett.\ Math.\ Phys.\  {\bf 24} (1992) 49.
%%CITATION = LMPHD,24,49;%%
%\cite{Aschieri:1999sg}
%\bibitem{Aschieri:1999sg}

P.~Aschieri,
``Real forms of quantum orthogonal groups, $q$-Lorentz groups in any dimension,''
Lett.\ Math.\ Phys.\  {\bf 49} (1999) 1
[math.qa/9805120].
%%CITATION = MATH.QA 9805120;%%


\bibitem{Varilly}
J.~C.~Varilly,
``Quantum symmetry groups of noncommutative spheres,''
Commun.\ Math.\ Phys.\  {\bf 221} (2001) 511
[math.qa/0102065].
%%CITATION = MATH-QA 0102065;%%


\bibitem{Woronowicz}
%\cite{Woronowicz:1989rt}
S.~L.~Woronowicz,
``Differential Calculus On Compact Matrix Pseudogroups (Quantum Groups),''
Commun.\ Math.\ Phys.\  {\bf 122} (1989) 125.
%%CITATION = CMPHA,122,125;%%


\end{thebibliography}
\end{document}